\documentclass[reqno,10pt,centertags]{amsart}
\usepackage{amsmath,amsthm,amscd,amssymb,latexsym,upref}

\newcommand{\bbN}{{\mathbb{N}}}
\newcommand{\bbR}{{\mathbb{R}}}

\newcommand{\bbC}{{\mathbb{C}}}

\newcommand{\calB}{{\mathcal B}}

\newcommand{\calD}{{\mathcal D}}
\newcommand{\calE}{{\mathcal E}}
\newcommand{\calF}{{\mathcal F}}
\newcommand{\calH}{{\mathcal H}}

\newcommand{\calK}{{\mathcal K}}

\newcommand{\no}{\nonumber}
\newcommand{\lb}{\label}

\newcommand{\ol}{\overline}
\newcommand{\ti}{\tilde  }
\newcommand{\wti}{\widetilde  }

\newcommand{\spec}{\text{\rm{spec}}}

\newcommand{\ran}{\text{\rm{ran}}}
\newcommand{\ind}{\text{\rm{ind}}}
\newcommand{\dom}{\text{\rm{dom}}}

\newcommand{\bi}{\bibitem}
\newcommand{\hatt}{\widehat}
\newcommand{\beq}{\begin{equation}}
\newcommand{\eeq}{\end{equation}}
\newcommand{\ba}{\begin{align}}
\newcommand{\ea}{\end{align}}
\renewcommand{\Re}{\text{\rm Re}}
\renewcommand{\Im}{\text{\rm Im}}
\renewcommand{\ln}{\text{\rm ln}}

\DeclareMathOperator{\tr}{tr}

\DeclareMathOperator*{\nlim}{n-lim}
\DeclareMathOperator*{\slim}{s-lim}

\allowdisplaybreaks
\numberwithin{equation}{section}
\newtheorem{theorem}{Theorem}[section]
\newtheorem{lemma}[theorem]{Lemma}

\newtheorem{hypothesis}[theorem]{Hypothesis}
\theoremstyle{definition}
\newtheorem{definition}[theorem]{Definition}
\newtheorem{example}[theorem]{Example}

\theoremstyle{remark}
\newtheorem{remark}[theorem]{Remark}

\begin{document}
\title[The Spectral Shift Operator]{The Spectral Shift 
Operator}
\author[Gesztesy, Makarov, and
Naboko]{Fritz Gesztesy,
Konstantin A.~Makarov, \\ and Serguei N.~Naboko }
\address{Department of Mathematics,
University of
Missouri, Columbia, MO
65211, USA}
\email{fritz@math.missouri.edu\newline
\indent{\it URL:}
http://www.math.missouri.edu/people/fgesztesy.html}
\address{Department of Mathematics, University of
Missouri, Columbia, MO
65211, USA}
\email{makarov@azure.math.missouri.edu}
\address{Department of Mathematical Physics, St. Petersburg 
State University, 198904 St. Petersburg,
Russia}
\email{naboko@snoopy.phys.spbu.ru}

%%%%%%%%%%%%%%%%%%%%%%%%%%%%%%%%%%%%%%%%%%%%%%%%%%%%%%%%
\begin{abstract}
We introduce the concept of a spectral shift operator and 
use it to derive Krein's spectral shift function for pairs
of self-adjoint operators. Our principal tools are
operator-valued Herglotz functions and their logarithms.
Applications to Krein's trace formula and to the
Birman-Solomyak spectral averaging formula are discussed.

\end{abstract}

\maketitle
%%%%%%%%%%%%%%%%%%%%%%%%%%%%%%%%%%%%%%%%%%%%%%%%%%%%%%%%
\section{Introduction}\lb{s1}
Krein's spectral shift function \cite{Kr62}, \cite{Kr83},
\cite{KY81} has received enormous attention in the past
due to its widespread applications in a variety of fields
including scattering theory (\cite{BW83}, Ch.~19, \cite{BK62},
\cite{BP98}, \cite{BY93}, \cite{BY93a}, \cite{JK78},
\cite{Ka78}, \cite{Kr83}, \cite{So93}, \cite{Ya92}, Ch.~8),
relative index theory (\cite{AS94}, \cite{BGGSS87}, 
\cite{BMS88},
\cite{GS88}, \cite{Mu98}, Ch.~X), spectral averaging
(\cite{Al87}, \cite{Bi94}, \cite{BS75}, \cite{DSS94},
\cite{Ja66}, \cite{Ja71}, \cite{Po98}, 
\cite{Si95}--\cite{SW86}) and its application to localization
properties of random Hamiltonians (\cite{BS98}, \cite{Ca83},
\cite{Ca84}, \cite{CL90}, Ch.~VIII, \cite{CH94},
\cite{CHM96}, \cite{KS98}--\cite{KS87},
\cite{PF92}, Ch.~V, \cite{Si85}, \cite{Si95}), eigenvalue
counting functions and spectral asymptotics (\cite{BP98},
\cite{Ch98}, \cite{EP97}, \cite{GKS95}, 
\cite{Pu97}--\cite{Pu98a}), semi-classical approximation 
(\cite{Ro92}--\cite{Ro99}, \cite{RS97}) 
and trace formulas for
one-dimensional Schr\"odinger and Jacobi operators
(\cite{GH95}--\cite{GHSZ95}, \cite{GKT96},
\cite{GS96}). Detailed reviews on the spectral shift
function and its applications up to 1993 were published by
Birman and Yafaev \cite{BY93}, \cite{BY93a}. These two papers 
contain extensive bibliographies and the interested reader 
will find many additional references therein. Historically, 
the concept of a spectral shift function was introduced by 
I.~M.~Lifshits \cite{Li52}, \cite{Li56}.

Our main contribution to this circle of ideas is the
introduction of a spectral shift operator
$\Xi(\lambda,H_0,H)$ for a.e. $\lambda\in\bbR,$ associated
with a pair of self-adjoint operators $H_0,$ $H=H_0+V$ with
$V\in\calB_1(\calH)$ ($\calH$ a complex separable Hilbert
space). In the special cases of sign-definite perturbations
$V\geq 0$ and $V\leq 0,$ $\Xi(\lambda,H_0,H)$ turns out to
be a trace class operator in $\calH,$ whose trace coincides
with Krein's spectral shift function $\xi(\lambda,H_0,H)$
for the pair $(H_0,H).$ While the special case $V\geq 0$ has
previously been studied by Carey \cite{Ca76}, our aim here
is to treat the case of general interactions $V$ by separately
introducing the positive and negative parts
$V_\pm=(|V|\pm V)/2$ of $V.$
In general, if $V$ is not sign-definite, then
$\Xi(\lambda,H_0,H)$ is not necessarily of trace class.
However, we will introduce trace class operators
$\Xi_\pm(\lambda)$ naturally associated with $V_\pm,$ 
acting in
distinct Hilbert spaces $\calH_{\pm},$ such that
\begin{equation}
\xi(\lambda,H_0,H)=\tr_{\calH_+}(\Xi_+(\lambda))
- \tr_{\calH_-}(\Xi_-(\lambda)) \text{ for a.e. }
\lambda\in\bbR. \lb{1.1}
\end{equation}
(An alternative approach to $\xi(\lambda,H_0,H),$ which
does not rely on separately introducing $V_+$ and $V_-,$
will be discussed elsewhere \cite{GM99}.)

Our main techniques are based on operator-valued Herglotz
functions (continuing some of our recent investigations
in this area \cite{GKMT98}, \cite{GMT98}, \cite{GT97})
and especially, on a detailed study of logarithms of Herglotz
operators in Section~\ref{s2}. In Section~\ref{s3} we 
introduce the
spectral shift operator $\Xi(\lambda,H_0,H)$ associated with
the pair $(H_0,H)$ and relate it to Krein's spectral
shift function $\xi(\lambda,H_0,H)$ and his celebrated
trace formula \cite{Kr62}. Finally, Section~\ref{s4} provides
an application to spectral averaging originally due to
Birman and Solomyak \cite{BS75} and hints at operator-valued
generalizations thereof. A number of additional
applications of this formalism will appear elsewhere
\cite{GM99}.

%%%%%%%%%%%%%%%%%%%%%%%%%%%%%%%%%%%%%%%%%%%%%%%%%%%%%%%
\section{Logarithms of Operator-Valued Herglotz Functions} 
\lb{s2}
%%%%%%%%%%%%%%%%%%%%%%%%%%%%%%%%%%%%%%%%%%%%%%%%%%%%%%%
The principal purpose of this section is to study 
operator-valued
 Herglotz functions and particularly their logarithms and
associated representation theorems.
In this  manner we shall obtain operator-valued
generalizations of some of the classical results on 
exponential Herglotz
representations studied by Aronszajn and Donoghue
\cite{AD56}.

In the following $\calH$ denotes a complex separable Hilbert 
space
with scalar product $(\,\cdot, \, \cdot)_{\calH}$
(linear in the second factor) and norm
$\vert\vert \cdot \vert\vert_\calH,$  $I_{\calH}$ the identity
operator in $\calH$,
$\calB(\calH)$ the Banach space of bounded linear
operators defined on $\calH$, $\calB_p(\calH), \,\, p\ge 1$
the standard Schatten-von Neumann ideals of
 $\calB(\calH)$ (cf., e.g., \cite{GK69}, \cite{Si79}) and
$\bbC_+ $ (resp., $\bbC_-$) the open complex upper
(resp., lower) half-plane. Moreover, real and
imaginary parts of a bounded operator $T\in \calB(\calH)$ are
defined as  usual by
$\Re (T)=(T+T^*)/2$,
$\Im (T)=(T-T^*)/(2i).$

%%%%%%%%%%%%%%%%%%%%%%%%%%%%%%%%%%%%%%%%%%%%%%%%%%%%%%%%%%%%%%
\begin{definition} \lb{d2.1}
$M:\bbC_+\to \calB(\calH)$ is called an
 {\it operator-valued  Herglotz function}
 (in short,
a {\it Herglotz operator}) if $M$ is analytic on $\bbC_+$ and
$\Im (M(z))\ge 0$ for all
$z\in \bbC_+$.
\end{definition}
%%%%%%%%%%%%%%%%%%%%%%%%%%%%%%%%%%%%%%%%%%%%%%%%%%%%%%%%%%%%%
\begin{theorem}\lb{t2.2} \mbox{\rm (Birman and Entina 
\cite{BE67},
de Branges \cite{de62}, Naboko \cite{Na87}--\cite{Na90}.)}
Let \linebreak $M:\bbC_+\to \calB(\calH)$ be a Herglotz 
operator.

\noindent(i) Then there exist bounded self-adjoint operators
 $A=A^*\in \calB(\calH),$
 $0\le B\in \calB(\calH)$, a Hilbert space
$ \calK\supseteq \calH,$ a self-adjoint
 operator $L=L^*$ in
$ \calK$, a bounded nonnegative operator
$0\le R\in \calB(\calK)$ with $R\vert_{\calK
\ominus \calH}=0$ such that
\begin{subequations}\lb{2.1}
\begin{align}
M(z)&=A+Bz+R^{1/2} (I_{\calK}+zL)(L-z)^{-1}R^{1/2}\vert_{\calH} 
\label{2.1a} \\
&= A+(B+R\vert_{\calH})z+(1+z^2)R^{1/2}(L-z)^{-1}
R^{1/2}\vert_{\calH}.
 \label{2.1b}
\end{align}
\end{subequations}

\noindent(ii) Let $p\ge1$. Then $M(z)\in \calB_p(\calH)$ 
for all
$z\in \bbC_+$ if and
only if $M(z_0)\in \calB_p(\calH)$ for some $z_0\in \bbC_+.$
In this case necessarily $A,B,R\in \calB_p(\calH).$

\noindent(iii) Let $M(z) \in \calB_1(\calH)$ for some (and
 hence for all) $z\in \bbC_+$. Then $M(z)$
 has normal boundary values $M(\lambda+i0)$ for (Lebesgue)
a.e. $\lambda\in \bbR$ in
every $\calB_p(\calH)$-norm, $p>1$. Moreover, let
$\{E_L(\lambda)\}_{\lambda\in\bbR}$ be the family of 
orthogonal spectral
projections of $L$ in $\calK.$ Then 
$R^{1/2}E_L(\lambda)R^{1/2}$ is
$\calB_1(\calH)$-differentiable for a.e. $\lambda\in\bbR$ and
denoting the derivative by 
$d(R^{1/2}E_L(\lambda)R^{1/2})/d\lambda,$
$\Im(M(z))$ has normal boundary values $\Im(M(\lambda+i0))$ 
for
a.e. $\lambda\in\bbR$ in $\calB_1(\calH)$-norm given by
\begin{equation}
\lim_{\varepsilon\downarrow 0}\|\pi^{-1}\Im(M(\lambda+
i\varepsilon))
- d(R^{1/2}E_L(\lambda)
R^{1/2}|_\calH)/d\lambda\|_{\calB_1(\calH)}
=0 \text{ a.e. } \lb{2.2a}
\end{equation}
\end{theorem}
%%%%%%%%%%%%%%%%%%%%%%%%%%%%%%%%%%%%%%%%%%%%%%%%%%%%%%%%%%%%%
Actually, the normal boundary values $M(\lambda+i0)$ in
Theorem~\ref{t2.2} (iii) can be replaced by nontangential 
ones.
Since this distinction will play no role in the remainder 
of this
paper we omit the corresponding details.

Theorem~\ref{t2.2} (i) follows by considering the quadratic 
form
 $(\varphi, M(z) \varphi)_{\calH}$,
 resulting in a scalar Herglotz function, in combination
with Naimark's dilation theorem (cf., e.g., Theorem 1 in 
Appendix I
of \cite{Br71}).
 Details can be found in \cite{Na87}. For 
Theorem~\ref{t2.2} (ii),
 (iii) we refer to
\cite{Na89}, \cite{Na90}. In particular, $p$ cannot be 
chosen equal
to $1$ in Theorem~\ref{t2.2} (iii) (cf. \cite{Na89}). Moreover,
if $M(z_0)\in\calB_p(\calH)$ for some $z_0\in\bbC_+$ and 
some $p>1,$
then $M(z)$ need not even have boundary values $M(\lambda+i0)$
in the weak topology of $\calH$ for a.e. $\lambda\in\bbR.$ 
In fact, the
quadratic form $(f,M(\lambda+i\varepsilon)f)_\calH$ may 
converge for
$f$ in a fixed set $\calD\subset\calH$ (independent of 
$\lambda$)
to $(f,M(\lambda+i0)f)_\calH$ for a.e. $\lambda\in\bbR$ as
$\varepsilon\downarrow 0,$ however, $M(\lambda+i0)$ may be 
a densely
defined unbounded operator in $\calH$ for a.e. $\lambda\in\bbR$
(cf. \cite{Na89}).

Originally, the existence of normal limits $M(\lambda+i0)$ for
a.e. $\lambda\in\bbR$ in $\calB_2(\calH)$-norm, in the special
case $A=0,$ $B=-R|_\calH,$ assuming $M(z)\in\calB_1(\calH),$ 
was
proved by de Branges \cite{de62}
in 1962. (The more general case in \eqref{2.1} can easily be
reduced to this special case.) In his paper \cite{de62}, 
de Branges
also proved the existence of normal limits $\Im(M(\lambda+i0))$
for a.e. $\lambda\in\bbR$ in $\calB_1(\calH)$-norm and obtained
\eqref{2.2a}.
These results and their implications on stationary scattering 
theory
were subsequently studied in detail by Birman and Entina
\cite{BE64}, \cite{BE67}. (Textbook representations of this
material can also be found in \cite{BW83}, Ch.~3.)

Invoking the family of orthogonal spectral projections
$\{E_L(\lambda)\}_{\lambda\in \bbR}$ of $L$, \eqref{2.1}
then yields the  familiar representation
\begin{equation}\lb{2.2}
M(z)= A+Bz+
\int_{\bbR} (1+\lambda^2)
d(R^{1/2}E_L(\lambda)
 R^{1/2}\vert_{\calH})
  ((\lambda-z)^{-1}
-\lambda(1+\lambda^2)^{-1}) ,
\end{equation}
where for our purpose it suffices to interpret
the integral in \eqref{2.2} in the weak sense.
Further results on representations of the type \eqref{2.2} 
can be found in
\cite{Br71}, Sect.~I.4  and  \cite{Sh71}.

Before we continue our investigations  on Herglotz operators 
a few
comments concerning our terminology
 are perhaps in order. The representation \eqref{2.2} in the 
special
 case of scalar Herglotz functions
is due to Riesz and Herglotz, respectively, Nevanlinna. The 
former
 authors discussed
 the analog
of \eqref{2.2} in the open unit disk $D\subset \bbC$,
 whereas the latter studied \eqref{2.2} in $\bbC_+$. As a
 consequence, functions of the type \eqref{2.2} are frequently 
called
Herglotz or Nevanlinna functions. Moreover,
 researchers in the former Soviet Union coined the term 
$R$-functions.
Here we follow the traditional
terminology in mathematical physics
 which seems to prefer the notion of Herglotz functions.

Since we are interested in logarithms of Herglotz operators, 
questions of
 their invertibility naturally arise.
The following result clarifies the situation.
%%%%%%%%%%%%%%%%%%%%%%%%%%%%%%%%%%%%%%%%%%%%%%%%%%%%%%%%%%%%
\begin{lemma}\lb{l2.3}
Suppose $M$ is a Herglotz operator with values in 
$\calB(\calH)$.
 If $M(z_0)^{-1}\in \calB(\calH)$ for some $z_0\in \bbC_+$ 
then
$M(z)^{-1}\in \calB(\calH)$ for all $z\in \bbC_+$.
\end{lemma}
%%%%%%%%%%%%%%%%%%%%%%%%%%%%%%%%%%%%%%%%%%%%%%%%%%%%%%%%%%%%%
\begin{proof}
Suppose there is a $z_1\in \bbC_+$ and a sequence
$\{e_n\}_{n\in \bbN}\subset \calH$ such that
$\vert \vert e_n\vert\vert_{\calH}  =1, $ $ n\in \bbN$ but
$\lim_{n\to \infty} \|M(z_1)e_n\|_{\calH} = 0.$ Then
$
\lim_{n\to \infty}
(e_n, \Im (M(z_1))e_n)_{\calH}=
\lim_{n\to \infty}
\vert\vert
 (\Im(M(z_1)))^{1/2}e_n\vert\vert_{\calH}=0$
and
hence
\begin{equation}\lb{2.3}
\lim_{n\to \infty}
\|\Im (M(z_1))e_n\|_{\calH} = 0.
\end{equation}
By \eqref{2.1a},
\begin{align}
\Im (M(z))&=\Im (z)( B+R^{1/2}(I_{\calK}+L^2)
((L-\Re(z))^2+(\Im (z))^2)^{-1}R^{1/2}\vert_{\calH})
 \no \\
&= \Im (z) (B+R^{1/2}C(z)R^{1/2}\vert_{\calH}),
 \label{2.4}
\end{align}
where $C(z)=(I_{\calK}+L^2)((L-\Re(z))^2+(\Im (z))^2)^{-1}>0$
is invertible in $\calK$, $C(z)^{-1}\in\calB(\calK)$, $ z\in 
\bbC_+$.
 Thus, \eqref{2.3} implies $\lim_{n\to \infty}
((e_n, Be_n)_{\calH} + (R^{1/2}\vert_{\calH} e_n, C(z_1)
R^{1/2}\vert_{\calH}$ $e_n)_{\calH})= 0$
and hence $\lim_{n\to \infty} \|B^{1/2}e_n\|_{\calH} = 0$, and
$\lim_{n\to \infty}\|Be_n\|_{\calH}= 0,$
$\lim_{n\to \infty}$ $\|R^{1/2}\vert_{\calH}e_n\|_{\calH}= 0.$
Consequently,
\begin{equation}\lb{2.8}
M(z_1)e_n=Ae_n+z_1Be_n
+(R^{1/2}(I_{\calK}+z_1L)(L-z_1)^{-1})R^{1/2}\vert_{\calH}e_n
\end{equation}
with
$R^{1/2}(I_{\calK}+z_1L)(L-z_1)^{-1} \in \calB(\calK)$ yields
$\lim_{n\to \infty}\|Ae_n\|_{\calH}= 0.$
Applying \eqref{2.1a} again we infer
$\lim_{n\to \infty}\|M(z)e_n\|_{\calH}= 0  \text{ for all }
z\in \bbC_+$
contradicting our hypothesis $M(z_0)^{-1}\in \calB(\calH)$.
This argument shows two facts. First, by choosing
$e_n = e \in \ker(M(z_1))$
it yields $\ker(M(z)) = \{0\} \text{ for all } z \in \bbC_+,$
and second, it proves the boundedness of
$M(z)^{-1}:\ran(M(z)) \to \calH,$
 $z
\in \bbC_+.$ In particular, one infers that $\ran(M(z))
= \dom(M(z)^{-1})$
is a
closed subspace of $\calH.$ By exactly the same argument one 
derives
$\{0\} = \ker(M(z)^*) = (\ran(M(z))^{\bot}, \quad z \in 
\bbC_+$ and
thus, $M(z)^{-1} \in \calB(\calH) \text{ for all } z \in 
\bbC_+.$
\end{proof}
%%%%%%%%%%%%%%%%%%%%%%%%%%%%%%%%%%%%%%%%%%%%%%%%%%%%%%%%%%%%%

Lemma~\ref{l2.3} admits the following generalization in 
connection
with parametrices, that is, generalized inverses, familar
from treatments
of the Calkin algebra. In the following we denote by 
$\calF(\calH)$
the algebra of finite-rank operators on $\calH$ and recall 
that $T$
is called a left (resp., right) parametrix of $M(z_0)$ if
$(TM(z_0)-I_\calH)\in\calF(\calH)$ (resp.,
$(M(z_0)T-I_\calH)\in\calF(\calH)$).

%%%%%%%%%%%%%%%%%%%%%%%%%%%%%%%%%%%%%%%%%%%%%%%%%%%%%%%%%%%%%%%
\begin{lemma} \lb{l2.3a} Suppose $M$ is a Herglotz operator 
with
values in $\calB(\calH).$ If $M(z_0)$ has a left {\rm (}resp.,
right{\rm \,)}
parametrix for some $z_0\in\bbC_+,$ then $M(z)$ has a left 
and right
parametrix for all $z\in\bbC_+.$
\end{lemma}
%%%%%%%%%%%%%%%%%%%%%%%%%%%%%%%%%%%%%%%%%%%%%%%%%%%%%%%%%%%%%%%
\begin{proof}
 From the proof of Lemma~\ref{l2.3} one infers
\begin{subequations} \lb{2.13a}
\begin{align}
\ker(M(z_1))&=\ker(M(z_2)) \text{ for all } z_1, z_2 \in\bbC_+,
\lb{2.13aa} \\
\ker(M(z)^*)&=\ker(M(z))=(\ran(M(z)))^{\bot}, \quad z\in\bbC_+
\lb{2.13ab}
\end{align}
\end{subequations}
and hence $M(z)$ can be decomposed as
\begin{equation}
M(z)=\begin{pmatrix}\widetilde M(z) & 0 \\ 0 & 0  \end{pmatrix},
\quad \calH =\widetilde \calH \oplus \ker(M(i)), \lb{2.13c}
\end{equation}
where
\begin{equation}
\ker(\widetilde M(z))=\{0\}, \quad \ol{\ran(\widetilde M(z))}
=\widetilde \calH, \quad z\in\bbC_+. \lb{2.13d}
\end{equation}
 From  $TM(z_0)-I_\calH=F\in\calF(\calH)$ one concludes
$\ker(M(z_0))\subseteq \ker(I_\calH + F)$
and hence
\begin{equation}
\dim(\ker(M(z)))<\infty \text{ for all } z\in\bbC_+. 
\lb{2.13g}
\end{equation}
Denote by $\widetilde P$ the orthogonal projection onto
$\widetilde \calH,$ $\widetilde \calH = \widetilde P \calH$
then $\widetilde P TM(z_0)\widetilde P = \widetilde P T 
\widetilde P
{\widetilde M(z_0)} =\widetilde P (I_\calH + F)\widetilde P
=I_{\ti \calH} + \widetilde F.$
In order to prove that ${\widetilde M(z_0)}^{-1}
\in\calB({\widetilde \calH})$ we argue as follows: Suppose 
there is
a sequence $\{e_n\}_{n\in\bbN}\subset {\wti \calH},$
$\|e_n\|_{\ti \calH}=1$ such that
\begin{equation}
\lim_{n\to\infty} \|\widetilde M(z_0) e_n\|_{\ti \calH} =0.
\lb{2.13j}
\end{equation}
Hence $\lim_{n\to\infty} \|\wti P T \wti P 
\wti M(z_0) e_n\|_{\ti \calH}
=\lim_{n\to\infty} \|(I_{\ti \calH} + 
\wti F)e_n\|_{\ti \calH}
=0.$
Next, let ${\wti E}_{-1}$ be the spectral projection
onto the geometric eigenspace of $\wti F$ corresponding to
the eigenvalue $-1,$ and
\begin{equation}
{\wti E}_{-1} \dotplus {\wti E}_{-1}^{\bot} = I_{\ti \calH}
\lb{2.13l}
\end{equation}
(here $\dotplus$ denotes the direct sum but not 
necessarily the
orthogonal direct sum $\oplus$). Then
\begin{equation}
\lim_{n\to \infty} \| (I_{\ti \calH} +
\wti F){\wti E}_{-1}^\bot e_n\|_{\ti \calH} =0. \lb{2.13m}
\end{equation}
Since $I_{\ti \calH} + \wti F$ is boundedly invertible on
${\wti E}_{-1}^\bot {\wti \calH},$ one infers
\begin{equation}
\lim_{n\to\infty} \|{\wti E}_{-1}^\bot e_n\|_{\ti \calH} =0.
\lb{2.13n}
\end{equation}
Introducing $\ti e_n = {\wti E}_{-1} e_n,$ $n\in\bbN,$ one
concludes from \eqref{2.13j}, \eqref{2.13l}--\eqref{2.13n}
that $\lim_{n\to\infty} \|\ti e_n\|_{\ti \calH} =1$,
$\lim_{n\to\infty} \|\wti M(z_0)\ti e_n\|_{\ti \calH}
=0.$ Since $\{\ti e_n\}_{n\in\bbN}$ is compact, there 
exists a
subsequence $\{\ti e_{n_k}\}_{k\in\bbN}\subseteq
\{\ti e_n\}_{n\in\bbN}$ such that
$\lim_{k\to\infty} \|\ti e_{n_k} - \ti e_0\|_{\ti \calH} =0$,
 $\ti e_0 \in \ran(\wti E_{-1})$, 
$\|\ti e_0\|_{\ti \calH} =1.$
Hence $\lim_{k\to\infty} \|\wti M(z_0) 
\ti e_{n_k}\|_{\ti \calH}=0$
implying the contradiction
$\wti M(z_0) \ti e_0 =0.$
Thus, $\wti M(z_0)^{-1} \in\calB(\wti \calH)$ and hence
$\wti M(z)^{-1} \in \calB(\wti \calH)$, $z\in\bbC_+$
by Lemma~\ref{l2.3}. In particular, $\ran(\wti M(z))=
\wti \calH,$ $z\in\bbC_+$ and hence $M(z)=
\left(\begin{smallmatrix}
\ti M(z)& 0 \\ 0 & 0  \end{smallmatrix}\right)$ has left and
right parametrices taking into account \eqref{2.13g}.
\end{proof}
%%%%%%%%%%%%%%%%%%%%%%%%%%%%%%%%%%%%%%%%%%%%%%%%%%%%%%%%%%%%%

Concerning boundary values at the real axis we also recall

%%%%%%%%%%%%%%%%%%%%%%%%%%%%%%%%%%%%%%%%%%%%%%%%%%%%%%%%%%%%%
\begin{lemma} \lb{l2.3b} \mbox{\rm (Naboko \cite{Na93}.)} 
Suppose
$(M-I_\calH)$ is a Herglotz operator with values in 
$\calB_1(\calH).$
Then the boundary values $M(\lambda+i0)$ exist for a.e.
$\lambda\in\bbR$ in $\calB_p(\calH)$-norm, $p>1$ and
$M(\lambda+i0)$ is a Fredholm operator for a.e. 
$\lambda\in\bbR$
with index zero a.e.,
\begin{equation}
\ind(M(\lambda+i0)) =0 \text{ for a.e. } \lambda\in\bbR.
\lb{2.13s}
\end{equation}
Moreover,
\begin{equation}
\ker(M(\lambda+i0)) =\ker(M(i)) = (\ran(M(\lambda+i0)))^\bot
\text{ for a.e. } \lambda \in\bbR. \lb{2.13t}
\end{equation}
In addition, if $M(z_0)^{-1}\in\calB(\calH)$ for some 
{\rm (}and hence
for all{\rm \,)} $z_0\in\bbC_+,$ then
\begin{equation}
M(\lambda+i0)^{-1}\in\calB(\calH) \text{ for a.e. }
\lambda\in\bbR. \lb{2.13u}
\end{equation}
\end{lemma}
%%%%%%%%%%%%%%%%%%%%%%%%%%%%%%%%%%%%%%%%%%%%%%%%%%%%%%%%%%%%

Generalizations and a counter example to \eqref{2.13u} if
$(M-I_\calH)\in\calB_1(\calH)$ is replaced by
$(M-I_\calH)\in\calB_p(\calH),$ $p>1$ can be found in 
\cite{Na93}.

Next, let $T$ be a bounded {\it dissipative} operator, 
that is,
\begin{equation}\lb{2.14}
T\in \calB(\calH), \quad \Im (T) \ge 0.
\end{equation}
In order to define the logarithm of $T$
we use the integral representation
\begin{equation}
\log (z) =-i \int_0^\infty d \lambda \,
 ((z+i\lambda)^{-1}- (1+i\lambda)^{-1}), \quad
z \neq-iy, \, y\geq 0, \lb{2.15}
\end{equation}
with a cut along the negative imaginary $z$-axis.
We use the symbol $\log(\cdot)$ in \eqref{2.15} in order to 
distinguish it
from the integral representation
\begin{equation}\lb{2.15a}
\ln(z) = \int_{-\infty}^0 d\lambda \, ((\lambda - z)^{-1} -
\lambda(1+\lambda^2)^{-1}), \quad z\in
\bbC\backslash (-\infty,0]
\end{equation}
with a cut along the negative real axis. Both representations 
will be used
later and it is easily verified that $\log(\cdot)$ and 
$\ln(\cdot)$ coincide
for $z\in\bbC_+.$ In particular,
one verifies that \eqref{2.15} and
\eqref{2.15a} are Herglotz functions,
 that is, they are analytic in $\bbC_+$
 and
\begin{equation}\lb{2.16}
0 < \Im (\log (z)) , \, \, \Im (\ln (z)) \, < \pi, \quad 
z\in \bbC_+.
\end{equation}
%%%%%%%%%%%%%%%%%%%%%%%%%%%%%%%%%%%%%%%%%%%%%%%%%%%%
\begin{lemma}\lb{l2.4}
Suppose $T\in \calB(\calH)$ is dissipative
and $T^{-1}\in \calB(\calH)$. Define
\begin{equation}\lb{2.17}
\log (T)=-i\int_0^\infty  d \lambda \,
((T+i\lambda)^{-1}-(1 +i\lambda)^{-1}I_\calH)
\end{equation}
in the sense of a $\calB(\calH)$-norm convergent Riemann 
integral.
Then

\noindent (i) $\log (T)\in \calB(\calH).$

\noindent (ii) If
$T=zI_{\calH},$ $z \in \bbC_+,$
 then $\log(T) =\log (z) I_{\calH}.$

\noindent (iii) Suppose
$\{ P_n\}_{n\in \bbN} \subset \calB(\calH)$ is a
family of  orthogonal finite-rank projections
in $\calH$ with $\slim_{n \to \infty}P_n=I_{\calH}.$ Then
$$
\slim_{n \to \infty}((I_{\calH}-P_n)+P_nTP_n)= T
$$
and
\begin{align}
&\slim_{n \to \infty}\log ((I_{\calH}- P_n) + 
P_n (T+i\varepsilon)P_n)
\no \\
&
=
\slim_{n \to \infty}P_n(\log (P_n (T+
i\varepsilon)P_n\vert_{P_n\calH})P_n =
 \log
(T+i\varepsilon I_{H}), \quad \varepsilon >0.\no
\end{align}

\noindent (iv) $
\nlim_{\varepsilon \downarrow 0 }\log (T+
i\varepsilon I_{\calH}) = \log (T).$

\noindent (v) $e^{\log(T)}=T.$
\end{lemma}
%%%%%%%%%%%%%%%%%%%%%%%%%%%%%%%%%%%%%%%%%%%%%%%%%%%%%%
\begin{proof} (i) Clearly $\log (T)\in \calB(\calH)$ since
\begin{equation}\lb{2.18}
\vert\vert\log (T)\vert\vert\le\int_0^\delta  d\lambda
\,(\vert\vert
(T+i\lambda)^{-1}\vert\vert +1)
+\int_\delta^\infty d\lambda \,(\vert \vert T\vert\vert+1) 
\lambda^{-2}
\end{equation}
using $\vert
\vert (T +
i\lambda)^{-1}\vert\vert \le \lambda^{-1}, $ $\lambda >0$. 
Moreover, by
$\vert\vert (T +
i\lambda)^{-1}\vert\vert\le \vert\vert T^{-1}\vert\vert
(1-\vert\vert T^{-1}\vert\vert \lambda)^{-1}$ for
 $0<\lambda <\vert\vert T^{-1}\vert\vert^{-1},$
choosing $\delta=2^{-1}\vert\vert T^{-1}\vert\vert^{-1}$ is
sufficient to bound the first integral in \eqref{2.18}.

(ii) is obvious from \eqref{2.15}.

(iii) For any $f\in \calH, $  $\varepsilon >0,$
\begin{align}
&
\bigg \|
\int_0^\infty d\lambda \,
(((I_{\calH}-P_n)+P_n(T+i\varepsilon)P_n +i\lambda)^{-1}
-(T+i\varepsilon +i\lambda)^{-1})f
\bigg \|
\no \\
&
\le
2\int_0^\delta d\lambda \,
((\varepsilon +\lambda)^{-1})  \vert\vert f \vert\vert
\no \\
&
+\int_\delta^N d\lambda \,
\vert
\vert (( I_{\calH} -P_n)+P_nTP_n -T+i\varepsilon 
(P_n -I_{\calH}))
(T+i\varepsilon +i\lambda)^{-1} f \vert\vert \delta^{-1}
\no \\
&
+\int_N^\infty d\lambda \,
\lambda^{-2} \|
 ((I_{\calH}-P_n)+P_n(T+i\varepsilon )P_n -T-i\varepsilon 
I_\calH)\|  \| f \|.
\lb{2.20}
\end{align}
Taking $\delta >0$ sufficiently small and $N$ sufficiently 
large such that the
first
and third integrals in
\eqref{2.20} are sufficiently small uniformly with respect 
to $n\in \bbN$ for
fixed $f\in \calH$, it suffices to let $n\to \infty$
in the second integral in \eqref{2.20} for fixed $\delta $ 
and $N$.

(iv) One estimates,
\begin{align}
&
\| \log (T+i\varepsilon) -\log (T)
\|
=\bigg \| \int_0^\infty d \lambda \,
 ((T+i\varepsilon +i\lambda)^{-1}- (T+i\lambda)^{-1}) \bigg \|
\no \\
&
\le \int_0^\delta
d \lambda \,
(\vert\vert T^{-1}\vert\vert (1-(\varepsilon +\lambda ) 
\vert\vert
T^{-1}\vert\vert )^{-1}+
\vert\vert T^{-1}
\vert\vert (1-\lambda \vert\vert T^{-1}\vert\vert )^{-1})
\no \\
&
+\varepsilon \int_\delta^\infty d\lambda \, \lambda^{-2}   
\le 2\delta
(\vert\vert
T^{-1}\vert\vert^{-1}-(\varepsilon +\delta))^{-1}+\varepsilon 
\delta^{-1}
\lb{2.21}
\end{align}
for $\varepsilon +\delta <\vert\vert T^{-1}\vert\vert^{-1}$. 
Taking $\delta \le
2^{-1}\vert \vert T^{-1}\vert\vert^{-1}$,
$\varepsilon < \delta/2$ and $\varepsilon \downarrow 0$ then 
yields the result
since
$\delta >0$ can be chosen arbitrarily small.

(v) clearly holds for dissipative $n\times n$ matrices and  
hence
\begin{equation}\lb{2.22}
e^{\log (( I_{\calH}-P_n)+P_n(T+i\varepsilon)P_n)}
=
(I_{\calH}-P_n)+P_n(T+i\varepsilon)P_n, \quad \varepsilon >0
\end{equation}
upon decomposing $\calH=P_n\calH\oplus (I_{\calH}-P_n)\calH$,
where $P_n$ are orthogonal rank-$n$ projections. Since
$\slim_{n \to \infty} e^{A_n}= e^A$
whenever $A_n, A\in \calB(\calH)$ and $\slim_{n \to 
\infty}A_n = A$
(simply use
$e^B=\sum_{m=0}^\infty B^m/(m!),$ $B\in \calB(\calH)$ and
$
A_n^m-A^m=\sum_{\ell=0}^{m-1} A_n^\ell (A_n-A) A^{m-1-\ell}),
$
one infers from (iii) and \eqref{2.22}
$e^{\log(T+i\varepsilon)}=T+i\varepsilon$, $\varepsilon >0.$
Together with (iv) this yields (v).
\end{proof}
%%%%%%%%%%%%%%%%%%%%%%%%%%%%%%%%%%%%%%%%%%%%%%%%%%%%%%%%%%%%
\begin{lemma} \lb{l2.5} Suppose $T\in \calB(\calH)$ is 
dissipative
and $T^{-1}\in \calB(\calH)$. Let $L$ be the minimal 
self-adjoint dilation
of $T$ in the Hilbert space $\calK\supseteq \calH$.
Then
\begin{equation}\lb{n2.25}
\Im (\log (T))= \pi P_{\calH} E_L(( -\infty,0))\vert_{\calH},
\end{equation}
where $P_{\calH}$ is the orthogonal projection in $\calK$ 
onto $\calH$ and
$\{ E_L(\lambda)\}_{\lambda\in \bbR}$ is the family of 
orthogonal
 spectral projections of $L$ in $\calK$.
In particular,
\begin{equation}\lb{n2.26}
0\le \Im (\log (T))\le \pi I_{\calH}.
\end{equation}
\end{lemma}
%%%%%%%%%%%%%%%%%%%%%%%%%%%%%%%%%%%%%%%%%%%%%%%%%%%%%%%%%%%%
\begin{proof}
By Sz.-Nagy's dilation theorem (see the corresponding result
in \cite{SF70}, Ch.~III, Sect.~2, Theorem 2.1  for 
contractions), one infers
\begin{equation}\lb{n2.27}
(T+i\lambda)^{-1}=P_{\calH} (L+i\lambda)^{-1}\vert_{\calH},
 \quad \lambda >0,
\end{equation}
where $L$ is  the minimal  self-adjoint dilation of $T$  
in  $\calK$.
Then the existence of
  $T^{-1}\in \calB(\calH)$ and \eqref{n2.27} yield
\begin{equation}\lb{n2.28}
E_L(\{0\})=0.
\end{equation}
In order to prove \eqref{n2.28} one can argue as follows. 
Consider the
contraction $S=(T-i)(T+i)^{-1}.$
According to Theorem 3.2, Ch.~I, Sect.~3 in \cite{SF70},
every contraction on the
Hilbert space $\calH$ corresponds to a decomposition of 
$\calH$ into
an orthogonal sum $\calH=\calH_0\oplus \calH_1$ of two 
reducing subspaces of
$S$ such that the part of $S$ on $\calH_0$ is unitary, and 
the part of $S$
on $\calH_1$ is completely non-unitary. Moreover, this 
decomposition is unique.
 (We recall that a contraction is called completely 
non-unitary if
there are no non-zero subspaces reducing this contraction 
to a unitary
operator.) Since $T$ is invertible, the unitary part of 
the contraction $S$
does
not have the eigenvalue $-1$. Since the minimal unitary 
dilation of a
completely non-unitary contraction has absolutely continuous 
spectrum (see
\cite{SF60}), we conclude that the minimal unitary dilation 
of $S$ does not
have
the eigenvalue $-1$ and hence the kernel of the minimal 
self-adjoint dilation
$L$ of the dissipative operator $T$ is trivial, proving
\eqref{n2.28}. Next, \eqref{2.17} implies
\begin{align}
&\Im (\log (T))=
-\int_0^\infty d\lambda \,\Re ((T+i\lambda)^{-1}-
(1+i\lambda)^{-1}I_{\calH}  )
\no \\
&=
-\int_0^\varepsilon d\lambda \,\Re
((T+i\lambda)^{-1}-(1+i\lambda)^{-1}I_{\calH}  )
\no \\
&
\quad -P_{\calH}\int_\varepsilon^\infty d\lambda \,\Re
((L+i\lambda )^{-1}-(1+i\lambda)^{-1}I_{\calK} ) 
\vert_{\calH}
\no \\
&=-\int_0^\varepsilon d\lambda \,\Re
((T+i\lambda)^{-1}-(1+i\lambda)^{-1}I_{\calH}  )
\no \\
&
\quad -P_{\calH}\int_\varepsilon^\infty  d\lambda \,
 (L(L^2+\lambda^2 )^{-1}-(1+\lambda^2)^{-1}I_{\calK} )
 \vert_{\calH}, \quad \varepsilon >0. \lb{n2.29}
\end{align}
 Let us prove that
\begin{equation}\lb{n2.30}
\slim_{\varepsilon \downarrow 0}
\int_\varepsilon^\infty  d\lambda \,
 L(L^2+\lambda^2 )^{-1}= (\pi / 2)(P_{L,+}-P_{L,-}),
\end{equation}
where
$P_{L,+}=E_L((0, \infty))$ and
$P_{L,-}=E_L(( -\infty,0))$ are the
 spectral projections of $L$ corresponding to the half-lines
 $(0, \infty)$ and $(-\infty, 0)$. Then \eqref{n2.25} 
and hence
\eqref{n2.26} follow from \eqref{n2.29},
\eqref{n2.30} and the fact that
\begin{equation}\lb{n2.31}
\lim_{\varepsilon \downarrow 0}
\bigg \|
\int_0^\varepsilon d\lambda \,
\Re ((T+i\lambda)^{-1}-(1+i\lambda)^{-1}I_{\calH}  )
\bigg \|
=0.
\end{equation}
Using the estimate
$| \int_\varepsilon^\infty
d\lambda \,\mu(\mu^2+\lambda^2)^{-1} | \le (\pi / 2),$ 
$\mu\in\bbR,$ $\varepsilon > 0$ by the spectral theorem
for the self-adjoint operator $L$,
we infer that the family of operators 
$\int_\varepsilon^\infty
 d\lambda \,
 L(L^2+\lambda^2 )^{-1}$ is uniformly bounded in 
$\varepsilon >0$ and
therefore,
it suffices to check the convergence \eqref{n2.30} on a 
dense set in
$\calK$. A natural candidate for this set is
 $\calD=  \{ \bigcup_{\delta >0} E_L(\bbR \backslash 
(-\delta, \delta))f
\, |
\, f \in \calK \},$
 which is dense in $\calK$ since by \eqref{n2.28} the kernel 
of  $L$
is trivial. For  $f\in \calD$ we have
\begin{align}
&\int_\varepsilon^\infty  d\lambda \,
 L(L^2+\lambda^2 )^{-1} f
=
\int_\varepsilon^\infty  d\lambda \,
\int_{\bbR\backslash  (-\delta, \delta)}
\mu(\mu^2+\lambda^2 )^{-1}dE_L(\mu) f
\no \\
&=
\int_\delta^\infty ((\pi/2)-
\arctan (\varepsilon/\mu))dE_L(\mu) f
-\int_{-\infty}^{-\delta} ((\pi/2)-
\arctan (\varepsilon/\vert
\mu\vert)) dE_L(\mu) f
\no \\
&
=(\pi / 2) (E_L((\delta, \infty))-E_L(( -\infty, 
-\delta)))f
\no \\
&
-
\int_\delta^\infty \arctan (\varepsilon/\mu ) dE_L(\mu) f
+
\int^{-\delta}_{-\infty}
 \arctan (\varepsilon/\vert\mu\vert ) dE_L(\mu) f \lb{n2.32}
\end{align}
for all  $\delta >0 $,  $\delta =\delta(f) $ small enough. 
For fixed $f\in
\calD$,
 going to the limit $\varepsilon \to 0$ in \eqref{n2.32},
we get
\begin{align}
\lim_{\varepsilon \downarrow 0}
\int_\varepsilon^\infty  d\lambda \,
 L(L^2+\lambda^2 )^{-1} f &=
(\pi / 2)(E_L((\delta, \infty))-E_L(( -\infty, \delta)))f
\no \\
&=
(\pi / 2)(P_{L,+}-P_{L,-})f,
\end{align}
 proving \eqref{n2.30}.
\end{proof}
%%%%%%%%%%%%%%%%%%%%%%%%%%%%%%%%%%%%%%%%%%%%%%%%%%%%%%%%%%%

Combining Lemmas \ref{l2.3} and \ref{l2.5}  one can 
prove the
following result.
%%%%%%%%%%%%%%%%%%%%%%%%%%%%%%%%%%%%%%%%%%%%%%%%%%
\begin{lemma}\lb{l2.6}
Suppose $M:\bbC_+\longrightarrow \calB(\calH)$ is
a Herglotz operator and $M(z_0)^{-1}\in \calB(\calH)$ 
for some
{\rm (}and hence for all{\rm \,)} $z_0\in \bbC_+$. Then
$\log (M):\bbC_+ \longrightarrow \calB(\calH)$ is a Herglotz 
operator
and
\begin{equation}\lb{2.32}
0\le \Im (\log (M(z)))\le \pi I_\calH, \quad z\in \bbC_+.
\end{equation}
\end{lemma}
\begin{proof}
Clearly
\begin{equation}\lb{2.31}
\log (M(z))=-i\int_0^\infty d\lambda \, ((M(z)+i\lambda)^{-1}
-(1+i\lambda)^{-1}I_\calH ), \quad z\in \bbC_+
\end{equation}
is analytic for $z\in \bbC_+$ since $M(z)^{-1}\in 
\calB(\calH)$
 for all $z\in \bbC_+$ by Lemma \ref{l2.3}. An application of
Lemma \ref{l2.5} then yields \eqref{2.32}.
\end{proof}
%%%%%%%%%%%%%%%%%%%%%%%%%%%%%%%%%%%%%%%%%%%%%%%%%%%%%%%%%%%%%%

Thus applying \eqref{2.1a} to $\log(M(z))$ one infers
\begin{equation}\lb{2.33}
\log(M(z))=C+Dz+\widetilde R^{1/2} (I_{\widetilde \calK}+
z \widetilde L) (\widetilde L-z)^{-1}
\widetilde R^{1/2}\vert_\calH,
\quad z\in
\bbC_+
\end{equation}
for some Hilbert space $\widetilde \calK \supseteq \calH,$ 
some bounded
 self-adjoint operators $C$, $0\le D\in \calB(\calH),$
a bounded nonnegative operator $0\le \widetilde R\in 
\calB(\calH)$
with  $\widetilde R\vert_{\widetilde{\calK}\ominus \calH } =0$,
and a self-adjoint operator
$\widetilde L=\widetilde L^*$ in
$\widetilde \calK$. By comparison with scalar Herglotz
functions one would expect that $D=0$. That this is indeed the
case is proved
next.
%%%%%%%%%%%%%%%%%%%%%%%%%%%%%%%%%%%%%%%%%%%%%%%%%%%%%%%%%%
\begin{lemma}\lb{l2.7}
$D=0$ in the representation \eqref{2.33} for
$\log(M(z)),$ $ z\in \bbC_+.$
\end{lemma}
\begin{proof}
Consider $z=iy, $ $y \uparrow \infty.$ Then
\begin{equation}\lb{2.34}
\vert
\vert y^{-1} \log(M(iy))-Di\vert\vert =O(y^{-1})
\end{equation}
and by \eqref{2.32},
\begin{align}
\pi y^{-1} I_\calH
&\ge y^{-1}\Im(\log (M(iy)))
\no \\
&=D
+y^{-1} \Im (\widetilde R^{1/2} (I_{\widetilde \calK}
+iy\widetilde L)(\widetilde L-iy)^{-1}
\widetilde R^{1/2}\vert_{\calH})
\no \\
&
=D+y^{-1}(... \ge 0 ...)\ge D\ge 0 \lb{2.35}
\end{align}
and hence $D=0$ since $y^{-1}$ can be chosen 
arbitrarily small.
\end{proof}

Introducing the family of orthogonal spectral projections
$\{E_{\widetilde L}(\lambda)\}_{\lambda\in \bbR}$ of 
$\widetilde L$ in $
\widetilde \calK$
one can thus rewrite \eqref{2.33} as
\begin{align}
&\log (M(z))=C
\no \\
&+\int_{\bbR} (1+\lambda^2)
d(\widetilde R^{1/2} E_{\widetilde L}(\lambda) 
\widetilde R^{1/2}\vert_{\calH})
 ((\lambda - z)^{-1}-\lambda(1+\lambda^2)^{-1}),
\quad z\in \bbC_+ \lb{2.36}
\end{align}
and hence
\begin{equation}\lb{2.37}
\Im (\log (M(z)))=\Im (z) \int_\bbR
d(\widetilde R^{1/2}
 E_{\widetilde L}(\lambda) \widetilde R^{1/2}\vert_{\calH})
(1+\lambda^2)\vert \lambda-z\vert^{-2}
\end{equation}
for $z\in \bbC_+,$ interpreting both integrals in 
\eqref{2.36} and
\eqref{2.37} in the weak
sense for
simplicity. Again by comparison with scalar Herglotz 
functions one
expects that
$d(\widetilde R^{1/2}
E_{\widetilde L}(\lambda)$ $ \widetilde R^{1/2}\vert_{\calH})$ 
is a
 purely absolutely continuous operator-valued measure on 
$\bbR$.
This is confirmed by the following result
($ \tr_\calH(\cdot)$  denotes the trace of trace class 
operators in $\calH$).
%%%%%%%%%%%%%%%%%%%%%%%%%%%%%%%%%%%%%%%%%%%%%%%%%%%%%%%%%%%
\begin{theorem}\lb{t2.8}
Suppose $M:\bbC_+ \longrightarrow \calB(\calH)$ is a Herglotz 
operator
and $M(z_0)^{-1} \in \calB(\calH)$ for some {\rm (}and hence
for all{\rm  \,)} $z_0\in\bbC_+.$ Then
there exists a family of bounded self-adjoint
weakly {\rm (}Lebesgue{\rm \,)} measurable operators
 $\{\Xi(\lambda) \}_{\lambda\in \bbR}\subset \calB(\calH),$
\begin{equation}\lb{2.38}
0\le \Xi(\lambda)\le I_\calH \text{ for a.e. }  
\lambda\in \bbR
\end{equation}
such that
\begin{equation}\lb{2.39}
\log(M(z))=C+
\int_\bbR d \lambda \, \Xi(\lambda)
((\lambda-z)^{-1}-\lambda(1+\lambda^2)^{-1}),
 \quad z\in \bbC_+
\end{equation}
the integral taken in the weak sense, where $C=C^* \in 
\calB(\calH).$
 Moreover, if $\Im (\log$ $(M(z_0)))\in \calB_1(\calH)$ for
 some {\rm (}and hence for all{\rm \,)}
 $z_0\in \bbC_+$, then
\begin{align}
&
0\le \Xi(\lambda)\in \calB_1(\calH)
 \text{ for a.e. }  \lambda\in \bbR,
\lb{2.40} \\
&
0\le \tr_{\calH}(\Xi (\cdot))\in 
L_{\rm{loc}}^1(\bbR;d\lambda), \quad
\int_\bbR d\lambda \, (1+\lambda^2)^{-1}
\tr_\calH (\Xi(\lambda))<\infty,
\lb{2.42}
\end{align}
and
\begin{equation}\lb{2.43}
\tr_\calH (\Im (\log (M(z))))=\Im (z) \int_\bbR
d\lambda \,
\tr_{\calH} (\Xi(\lambda))
\vert \lambda-z\vert^{-2} , \quad z\in \bbC_+.
\end{equation}
\end{theorem}
%%%%%%%%%%%%%%%%%%%%%%%%%%%%%%%%%%%%%%%%%%%%%%%%%%%%%%%%%%%
\begin{proof}
Let $f\in \calH$ and denote
\begin{equation}
d\omega_f(\lambda)=
(1+\lambda^2) d(\widetilde R^{1/2}f, E_{\widetilde L}(\lambda)
 \widetilde R^{1/2}f)_\calH=
(1+\lambda^2) d
\vert\vert E_{\widetilde L}(\lambda)
\widetilde R^{1/2}f \vert\vert^2_\calH.
\lb{2.44}
\end{equation}
Then
\begin{equation}\lb{2.45}
0\le (f, \Im (\log (M(z)))f)_\calH\le \pi \vert\vert 
f \vert\vert^2_\calH,
\quad
z\in \bbC_+, \quad f\in \calH
\end{equation}
proves that $d\omega_f$ is purely absolutely continuous,
$d\omega_f=d\omega_{f,ac}$ by standard arguments (see, 
e.g., \cite{AD56}).
 Thus,
\begin{equation}
d\omega_f(\lambda)=\xi_f(\lambda)d \lambda  \text{ with }
 0\le \xi_f (\lambda ) \le \vert\vert f 
\vert\vert^2_\calH\quad
 \text{ for all }
  f\in \calH  \text{ and a.e. } \lambda\in \bbR.\lb{2.46}
\end{equation}
By \eqref{2.44}, $\xi_f(\lambda)$ defines a quadratic form
\begin{equation}\lb{2.47}
\xi_f(\lambda)=(f, \Xi(\lambda) f)_\calH
\text{ for some }
0\le \Xi(\lambda ) \le I_\calH
\text{ and a.e. }
 \lambda\in \bbR
\end{equation}
proving \eqref{2.38} and \eqref{2.39}. The representation 
\eqref{2.33}
(with $D=0$) implies
\begin{equation}
\Im (\log (M(z))=
\widetilde R^{1/2}
(I_\calK+\widetilde L^2)((\widetilde L-\Re(z))^2+
(\Im (z))^2)^{-1}
\widetilde R^{1/2}\vert_\calH \lb{2.48}
\end{equation}
and hence $\Im (\log (M(z_0)))\in \calB_1(\calH)$ for 
some $z_0 \in\bbC_+$
 implies
 $\Im (\log (M(z)))\in \calB_1(\calH)$ for all 
$z\in \bbC_+.$
In particular,
\begin{equation}\lb{2.49}
\Im (\log (M(i)))=\int_\bbR d\lambda \,(1+\lambda^2)^{-1} 
\Xi(\lambda)
\end{equation}
then proves
\begin{equation}\lb{2.50}
\tr_\calH (\Im (\log (M(i))))=
\int_\bbR  d\lambda \, (1+\lambda^2)^{-1}
\tr_\calH (\Xi(\lambda))<\infty
\end{equation}
by the monotone convergence  theorem. Hence
 $0\le \tr_\calH(\Xi(\lambda))=
\vert\vert \Xi(\lambda)\vert\vert_1<\infty$
 for a.e. $\lambda\in \bbR$ completes the proof.
\end{proof}
%%%%%%%%%%%%%%%%%%%%%%%%%%%%%%%%%%%%%%%%%%%%%%%%%%%%%%%%
\begin{remark}\lb{r2.9}
For simplicity we focused on dissipative operators thus far.
Later  we will
also encounter
operators  $S\in \calB(\calH)$ with $-S$ dissipative, 
that is,
$\Im (S) \le 0$.
In this case $S^*$ is dissipative and one can simply define 
$\log(S)$ by
\begin{equation}
\log(S) = (\log(S^*))^*,   \lb{2.55a}
\end{equation}
with $\log(S^*)$ defined as in \eqref{2.17}.
 Moreover,
\begin{align}
&\log (\hatt M(z))=\hatt C -
\int_\bbR d \lambda \, \hatt \Xi(\lambda)
((\lambda-z)^{-1}-\lambda(1+\lambda^2)^{-1}), \quad
z\in \bbC_+, \lb{2.52} \\
&\widehat C =\widehat C^* \in \calB(\calH)
\text{ and } 0 \le \hatt \Xi(\lambda )\le I_\calH  
\text{ for a.e. }
\lambda\in \bbR, \lb{2.52a}
\end{align}
whenever $\hatt M$ is analytic in $\bbC_+$ and
 $\Im (\hatt M(z))\le 0,$ $z\in \bbC_+.$
\end{remark}
%%%%%%%%%%%%%%%%%%%%%%%%%%%%%%%%%%%%%%%%%%%%%%%%
\begin{remark}\lb{r2.10}
Theorem \ref{t2.8} represents the operator-valued 
generalization of the
exponential Herglotz representation for scalar Herglotz 
functions studied in
detail
by Aronszajn and Donoghue
\cite{AD56} (see also Carey and Pepe
\cite{CP73}). Theorem
\ref{t2.8} is not the
 first attempt in this direction. Carey \cite{Ca76}, 
in 1976,
considered the case $M(z)=I_\calH +K^*(H_0-z)^{-1}K$ (i.e.,
$A-R^{1/2} L  R^{1/2}\vert_\calH  =I_\calH,$
 $B=0,$  $(1+L^2)^{1/2} R^{1/2}\vert_\calH=K$ when 
compared to
\eqref{2.1a}) and established
\begin{equation}\lb{2.56}
M(z)=\exp \bigg (\int_\bbR d\lambda \,B(\lambda) 
(\lambda-z)^{-1}
\bigg )
\end{equation}
for  a summable operator function $B(\lambda)$, 
$0\le B(\lambda)\le I_\calH$
(i.e., the analog of $\Xi(\lambda)$ in \eqref{2.39}). 
Although Carey's proof
also uses Naimark's dilation theorem as described
in Theorem 1 of Appendix I of \cite{Br71},
it is different from ours and
does not utilize the integral representation \eqref{2.17} 
for logarithms.
\end{remark}
%%%%%%%%%%%%%%%%%%%%%%%%%%%%%%%%%%%%%%%%%%%%%%%%%%%%%%%%%%%%
\begin{remark} \lb{r2.11}
At first glance it may seem that we have been a bit pedantic 
in introducing
various branches of the logarithm in \eqref{2.15} and 
\eqref{2.15a}. Actually,
these branches are just a special case of the following 
family of branches
\begin{equation}
\log_\alpha(z) =
\int_{\gamma_\alpha} d\zeta\,((z-\zeta)^{-1} -(1-\zeta)^{-1}),
\quad  z\in\bbC\backslash\gamma_\alpha, \lb{2.60a}
\end{equation}
where $\gamma_\alpha$ denotes the ray,
 $\gamma_\alpha = \{\zeta\in\bbC_+\,|\, \zeta = 
re^{i\alpha}, \, 0\leq r <
\infty, \, \alpha \in [\pi,2\pi] \}.$
In particular, $\log(\cdot)=\log_{3\pi/2}(\cdot)$ and since
$\int_{-\infty}^0
d\lambda\,
((1 - \lambda)^{-1} - \lambda(1 + \lambda^2)^{-1}) = 0,$
one infers
$\ln(\cdot)=\log_\pi(\cdot).$ That some care has to be taken 
in connection with
a consistent choice of branches especially for operator-valued 
branches of the
logarithm is illustrated in the following Remark~\ref{r2.12} 
and in
Remark~\ref{r3.2a}.
\end{remark}
\begin{remark}\lb{r2.12} In view of our applications in
Sections 3 and 4 it seems
worthwhile to recall in connection with our hypothesis
$\Im (\log(M(z_0))) \in
\calB_1(\calH)$
in Theorem \ref{t2.8} that if
$A\in\calB_1(\calH)$
 and
 $\text{det}_\calH(I_\calH + A)\neq 0,$
 then
 $\log(I_\calH +A)\in\calB_1(\calH)$
by \eqref{2.17} (using
$(I_\calH + A +i\lambda)^{-1} -
 (1 + i\lambda)^{-1}I_\calH
= (1 + \lambda)^{-1}(I_\calH + A +
i\lambda)^{-1}A.$) Moreover, $\text{det}_\calH (I_\calH + A)
= \prod_{n\in\bbN}(1 + \lambda_n (A)),$ where $\lambda_n(A)$
denote the eigenvalues of $A$ repeated according to their
algebraic multiplicity, then shows
\begin{equation}
\tr_\calH (\log(I_\calH + A)) = 
\log(\text{det}_ \calH (I_\calH + A)).
\lb{2.62}
\end{equation}
In fact, \eqref{2.62} holds for any branch 
$\log_\alpha (\cdot)$ and hence, in
particular, for the branch $\ln(\cdot)$ on either side
of \eqref{2.62}. (Here
$\text{det}_\calH (\cdot)$ denotes
the Fredholm determinant for operators in $\calH.$)
\end{remark}

%%%%%%%%%%%%%%%%%%%%%%%%%%%%%%%%%%%%%%%%%%%%%%%%%%%%%%%%%%%%%%
%%%%%%%%%%%%%%%%%%%%%%%%%%%%%%%%%%%%%%%%%%%%%%%%%%%%%%%%%%%%%%
\section{The Spectral Shift Operator}\lb{s3}

The main purpose of this section is to introduce the concept
of a spectral shift operator (cf.~Definition~\ref{d3.4a}) 
and a
new approach to Krein's basic trace formula \cite{Kr62}.

Let $\calH$ be a complex separable Hilbert space and assume 
the following
hypothesis for the remainder of this section.
%%%%%%%%%%%%%%%%%%%%%%%%%%%%%%%%%%%%%%%%%%%%%%%%%%%%%%%%%%%%%%%
\begin{hypothesis}\lb{h3.1}
Let  $H_0$ be a self-adjoint operator in $\calH$
with domain $\dom (H_0)$, $J$ a bounded self-adjoint operator
with $J^2=I_{\calH}$, and $K\in \calB_2(\calH)$ a 
Hilbert-Schmidt operator.
\end{hypothesis}
%%%%%%%%%%%%%%%%%%%%%%%%%%%%%%%%%%%%%%%%%%%%%%%%%%%%%%%%%%%%%%%

Introducing
\begin{equation}\lb{3.1}
V=KJK^*
\end{equation}
we define the self-adjoint operator
\begin{equation}\lb{3.2}
H=H_0+V, \quad \dom(H)=\dom(H_0)
\end{equation}
in $\calH$.

Given Hypothesis \ref{h3.1} we decompose $\calH$ and $J$ 
according to
\begin{equation}\lb{3.3}
J=\begin{pmatrix} I_+ & 0\\ 0&
-I_-\end{pmatrix}, \quad \calH=\calH_+\oplus \calH_-,
\end{equation}
\begin{equation}\lb{3.4}
J_+=\begin{pmatrix} I_+ & 0\\ 0&
0\end{pmatrix}, \quad
J_-=\begin{pmatrix} 0& 0\\ 0&
I_-\end{pmatrix}, \quad
J=J_+-J_-,
\end{equation}
with $I_\pm$ the identity operator in $\calH_\pm$. Moreover, 
we introduce
the following bounded operators
\begin{align}
\Phi (z)
&=J+K^*(H_0-z)^{-1}K:\calH \rightarrow \calH, \lb{3.5}\\
\Phi_+(z)
&=I_++J_+ K^*(H_0-z)^{-1}K\vert_{\calH_+}:
\calH_+ \rightarrow\calH_+, \lb{3.6} \\
\widetilde \Phi_-(z)
&=I_- - J_-K^*(H_+-z)^{-1}K\vert_{\calH_-}:
\calH_- \rightarrow \calH_-, \lb{3.7}
\end{align}
for $z\in \bbC\backslash \bbR,$ where
 \begin{equation}\lb{3.7'}
V_+=KJ_+K^*,
\end{equation}
\begin{equation}\lb{3.7''}
H_+=H_0+V_+, \quad \dom (H_+)=\dom (H_0).
\end{equation}
%%%%%%%%%%%%%%%%%%%%%%%%%%%%%%%%%%%%%%%%%%%%%%%%%%%%%%%%%%%%
\begin{lemma}\lb{l3.2}
Assume Hypothesis \ref{h3.1}. Then $\Phi$, $\Phi_+,$
 and $-\widetilde \Phi_-$ are Herglotz operators in
$\calH,$  $\calH_+,$ and
$\calH_-,$ respectively. In addition {\rm (}$z\in
\bbC\backslash \bbR${\rm )},
\begin{align}
\Phi (z)^{-1}
&=J-JK^*(H-z)^{-1}KJ, \lb{3.8} \\
\Phi_+ (z)^{-1}
&=I_+-J_+K^*(H_+-z)^{-1}K\vert_{\calH_+}, \lb{3.9} \\
\widetilde \Phi_- (z)^{-1}
&=I_- +J_-K^*(H-z)^{-1}K\vert_{\calH_-}.
\lb{3.10}
\end{align}
\end{lemma}
\begin{proof}
It suffices to consider $\Phi_+.$ Since  $I_+$ and $H_0$ are
self-adjoint,
$\Phi_+(z)=(J_++J_+K^*(H_0-z)^{-1}KJ_+)\vert_{\calH_+}$ is 
clearly
analytic in $\bbC\backslash\bbR$ and satisfies
$\Im (\Phi_+(z))\ge 0$ for $z\in \bbC_+.$
Relation \eqref{3.9} is then an elementary consequence of the
second resolvent equation,
\begin{subequations}\lb{3.12}
\begin{align}
(H_+-z)^{-1}&=
(H_0-z)^{-1}-(H_0-z)^{-1}V_+(H_+-z)^{-1}
 \lb{3.12a} \\
&=(H_0-z)^{-1}-(H_+-z)^{-1}V_+(H_0-z)^{-1},
\lb{3.12b}
\end{align}
\end{subequations}
the fact $J_+^2=J_+$, and simply follows by multiplying 
the right-hand
sides of
\eqref{3.6} and \eqref{3.9} in either order.
\end{proof}
%%%%%%%%%%%%%%%%%%%%%%%%%%%%%%%%%%%%%%%%%%%%%%%%%%%%%%%%%%
\begin{remark} \lb{r3.2a}
In the following we need to make use of the formula
\begin{equation}
d \tr_\calK (\log(I_\calK + F(z)))/dz = 
\tr_\calK (F'(z)(I_\calK +
F(z))^{-1}).  \lb{3.15a}
\end{equation}
This result is proven, for instance, in \cite{Di69}, 
Ch.~I, Sect.~6.11 or \cite{GK69},
Sect.~IV.1 for analytic $F(\cdot)\in\calB_1(\calK)$ 
in some region
$\Omega\subset\bbC$ such that $(I_\calK + 
F(\cdot))^{-1}\in\calB(\calK)$ in
$\Omega$ using the standard branch $\ln(\cdot).$ In 
this case equation
\eqref{3.15a} is first proven in the finite-dimensional 
case, followed by
taking the limit $n\to\infty$ upon replacing $F(z)$ by 
$P_nF(z)P_n,$ with $P_n$
a family of orthogonal projections
in $\calK$ strongly converging to $I_\calK$
as $n\to\infty$ (cf., \cite{GK69}, p.~163). This strategy 
of proof extends to
all branches $\log_\alpha(\cdot)$ introduced in 
Remark~\ref{r2.11}. More generally, we have the following 
result,
\begin{equation}
d \tr_\calK (\varphi( F(s)))/ds = 
\tr_\calK (\varphi'(F(s)) F'(s)),
  \lb{3.15b}
\end{equation}
with $F(s)\in \calB_1(\calK)$ defined on an interval 
$s_1\le s\le s_2$,
continuously differentiable in $\calB_1(\calK)$-norm, 
and $\varphi(z)$ any 
scalar function, holomorphic in some domain 
$\calD\subset \bbC$
with a Jordan boundary and
$\text{spec} (F(s)) \subset \calD$ for all 
$s\in [s_1, s_2].$
\end{remark}
%%%%%%%%%%%%%%%%%%%%%%%%%%%%%%%%%%%%%%%%%%%%%%%%%%%%%%%%%%%%

%%%%%%%%%%%%%%%%%%%%%%%%%%%%%%%%%%%%%%%%%%%%%%%%%%%%%%%%%%%%%
\begin{lemma}\lb{l3.3}
Assume Hypothesis \ref{h3.1} and $\bbC\backslash\bbR$. Then
\begin{subequations} \lb{3.13}
\begin{align}
\tr_\calH((H_0-z)^{-1}-(H_+-z)^{-1})
&=d \tr_{\calH_+} ( \log (\Phi_+(z)))/dz, \lb{3.13a} \\
\tr_\calH((H_+-z)^{-1}-(H-z)^{-1})
&=d \tr_{\calH_-} ( \log (\widetilde\Phi_-(z)))/dz. \lb{3.13b}
\end{align}
\end{subequations}
\end{lemma}
%%%%%%%%%%%%%%%%%%%%%%%%%%%%%%%%%%%%%%%%%%%%%%%%%%%%%%%%%%%%%%
\begin{proof}
It suffices to prove \eqref{3.13a}.
 Applying \eqref{3.12}  repeatedly one infers
\begin{align}
&\tr_\calH((H_0-z)^{-1}-(H_+-z)^{-1})
 \no \\
&=\tr_\calH((H_0-z)^{-1}V_+[(H_0-z)^{-1}-
(H_+-z)^{-1}V_+(H_0-z)^{-1}])
\no \\
&=\tr_\calH
(J_+K^*(H_0-z)^{-2}KJ_+[J_+-J_+K^*(H_+-z)^{-1}KJ_+])
\no\\
&=\tr_\calH
(J_+\Phi_+'(z)J_+^2\Phi_+(z)^{-1}J_+)=\tr_{\calH_+}
(\Phi_+'(z)\Phi_+(z)^{-1})
\no\\
&=d \tr_{\calH_+} (\log (\Phi_+(z)))/dz, \quad z \in 
\bbC\backslash \bbR,
\lb{3.15}
\end{align}
where we used $J_+^2=J_+$,
\begin{equation}\lb{3.15'}
\tr_{\calH}(AB)=\tr_\calH(BA)
\end{equation}
for $A,B \in \calB(\calH)$ with $AB, $ $BA\in \calB_1(\calH)$
(cf. Corollary 3.8 in \cite{Si79}), and
\eqref{3.15a}.
\end{proof}
%%%%%%%%%%%%%%%%%%%%%%%%%%%%%%%%%%%%%%%%%%%%%%%%%%%%%%%%%%%%%

Next, applying Theorem \ref{t2.8} and Remark \ref{r2.9} to
$\Phi_+(z)$ and $\widetilde \Phi_-(z)$ one infers the 
existence of two families
of bounded operators
$\{\Xi_\pm (\lambda)\}_{\lambda \in \bbR}$ defined for 
(Lebesgue) a.e.
$\lambda\in \bbR$ and
 satisfying
\begin{align}
&0\le \Xi_{\pm}(\lambda)\le I_\pm, \quad
\Xi_\pm (\lambda) \in \calB_1(\calH_\pm)
\text{ for a.e. } \lambda\in \bbR, \lb{3.16} \\
&
\vert\vert
\Xi_\pm(\cdot)\vert\vert_1\in L^1(\bbR;(1+
\lambda^2)^{-1} d\lambda)
\no
\end{align}
and
\begin{subequations} \lb{3.17}
\begin{align}
\log (\Phi_+(z))
&=
\log (I_++J_+K^*(H_0-z)^{-1}K\vert_{\calH_+})
\no\\
&=
C_+ +
\int_\bbR d\lambda \, \Xi_+(\lambda)
((\lambda-z)^{-1}-\lambda (1+\lambda^2)^{-1}),
\lb{3.17a} \\
\log (\widetilde \Phi_-(z))
&=
\log (I_--J_-K^*(H_+-z)^{-1}K\vert_{\calH_-})
\no \\
&=C_- -
\int_\bbR d\lambda  \,\Xi_-(\lambda)
((\lambda-z)^{-1}-\lambda (1+\lambda^2)^{-1})
\lb{3.17b}
\end{align}
\end{subequations}
for $z\in \bbC\backslash \bbR, $ with $C_\pm =C_\pm^*\in
\calB_1(\calH)$.

Equations \eqref{3.17} motivate the following
%%%%%%%%%%%%%%%%%%%%%%%%%%%%%%%%%%%%%%%%%%%%%%%%%%%%%%%%%%%
\begin{definition} \lb{d3.4a}
$\Xi_+(\lambda)$ (resp., $\Xi_-(\lambda)$) is called the
{\it spectral shift operator} associated with
$\Phi_+(z)$ (resp., $\widetilde \Phi_-(z)$). Alternatively, 
we will
refer to
$\Xi_+(\lambda)$ as the spectral shift operator
associated with the pair $(H_0,H_+)$ and occasionally
use the notation $\Xi_+(\lambda,H_0,H_+)$ to stress the 
dependence on
$(H_0,H_+),$ etc.
\end{definition}
%%%%%%%%%%%%%%%%%%%%%%%%%%%%%%%%%%%%%%%%%%%%%%%%%%%%%%%%%%%

Moreover, we introduce
\begin{equation}\lb{3.19}
\xi_\pm (\lambda) =\tr_{\calH_\pm}( \Xi_{\pm}(\lambda) ), \quad
0\le \xi_\pm\in L^1(\bbR;(1+\lambda^2)^{-1}d\lambda)
\text{ for a.e. } \lambda \in \bbR.
\end{equation}

Actually, taking into account the simple  behavior of
$\Phi_+(iy)$ and
$\widetilde\Phi_-(iy)$ as $
\vert y\vert \to \infty,$ one can improve
\eqref{3.17a} and \eqref{3.17b} as follows.
%%%%%%%%%%%%%%%%%%%%%%%%%%%%%%%%%%%%%%%%%%%%%%%%%%%%%%%%%%%
\begin{lemma}\lb{l3.4}
Assume Hypothesis \ref{3.1}
 and define $\xi_\pm$ as in \eqref{3.19}. Then
\begin{equation}\lb{3.19''}
0\le \xi_\pm \in L^{1}(\bbR; d\lambda),
\end{equation}
and \eqref{3.17a} and \eqref{3.17b} simplify to
\begin{subequations} \lb{3.25}
\begin{align}
\log(\Phi_+(z))
&=\int_\bbR d\lambda \, \Xi_+(\lambda)(\lambda -z)^{-1}, 
\lb{3.25a}
\\
\log(\widetilde\Phi_-(z))
&=-\int_\bbR d\lambda \, \Xi_-(\lambda)(\lambda -z)^{-1}.
\lb{3.25b}
\end{align}
\end{subequations}
Moreover, for a.e. $\lambda\in\bbR,$
\begin{subequations} \lb{3.26a}
\begin{align}
\lim_{\varepsilon\downarrow 0}\|\Xi_+(\lambda)-
\pi^{-1}\Im(\log(\Phi_+(\lambda+i\varepsilon)))\|
_{\calB_1(\calH_+)} &=0, \lb{3.26aa} \\
\lim_{\varepsilon\downarrow 0}\|\Xi_-(\lambda) +
\pi^{-1}\Im(\log(\wti \Phi_-(\lambda+i\varepsilon)))\|
_{\calB_1(\calH_-)} &=0. \lb{3.26ab}
\end{align}
\end{subequations}
\end{lemma}
%%%%%%%%%%%%%%%%%%%%%%%%%%%%%%%%%%%%%%%%%%%%%%%%%%%%%%%%%
\begin{proof}
It suffices to consider $\xi_+(\lambda)$ and $\Phi_+(z).$ 
Since
\begin{equation}\lb{3.19b}
\vert\vert
\log (\Phi_+(y) )\vert\vert_1=O
(\vert y
\vert^{-1}) \text{ as }
\vert y\vert \to  \infty
\end{equation}
by the Hilbert-Schmidt hypothesis on $K$ and the fact
$\vert\vert (H_0-iy)^{-1}\vert\vert =O
(\vert y
\vert^{-1}) $
as $\vert y\vert \to \infty$, the scalar Herglotz function
$\tr_{\calH_+}(\log (\Phi_+(z)))$ satisfies
\begin{equation}\lb{3.19c}
\vert\tr_{\calH_+}(\log (\Phi_+(z)))\vert
 =O
(\vert y
\vert^{-1})
\text{ as } \vert y\vert \to \infty.
\end{equation}
By standard results (see, e.g., \cite{AD56}, \cite{KK74}),
\eqref{3.19c} yields
\begin{equation}\lb{3.19d}
\tr_{\calH_+}(\log (\Phi_+(z)))=
\int_\bbR d \omega_+(\lambda) (\lambda-z)^{-1}, \quad
z\in \bbC\backslash\bbR,
\end{equation}
where $\omega_+$ is a finite measure,
\begin{equation}\lb{3.19e}
\int_\bbR d \omega_+(\lambda)=-i\lim_{y\uparrow\infty}
(y\tr_{\calH_+}(\log (\Phi_+(z))))<\infty.
\end{equation}
Moreover, the fact that $\Im(\log(\Phi_+(z)))$ is uniformly 
bounded
with respect to $z\in\bbC_+$ yields that $\omega_+$ is 
purely absolutely
continuous,
\begin{equation}\lb{3.19f}
d \omega_+(\lambda)=\xi_+(\lambda)d\lambda,
\quad \xi_+\in L^1(\bbR ; d\lambda),
\end{equation}
where
\begin{align}
\xi_+(\lambda)
&=\pi^{-1} \lim_{\varepsilon \downarrow 0}
(\Im (\tr_{\calH_+}(\log(\Phi_+(\lambda+i\varepsilon)))))
=\tr_{\calH_+}(\Xi_+(\lambda))
\no \\
&
=\pi^{-1}\lim_{\varepsilon \downarrow 0}
(\Im (\log (\text{det}_{\calH_+} (\Phi_+(\lambda+
i\varepsilon)))))
\text{ for a.e. }
\lambda\in \bbR. \lb{3.19g}
\end{align}
In order to prove \eqref{3.26aa} we first observe that
$\Im(\log(\Phi_+(\lambda+i\varepsilon)))$ takes on boundary 
values
$\Im(\log(\Phi_+(\lambda+i0)))$ for a.e. $\lambda\in\bbR$ in
$\calB_1(\calH_+)$-norm by \eqref{2.2a}. Next, choosing an 
orthonormal
system $\{e_n\}_{n\in\bbN}\subset\calH_+,$ we recall that
the quadratic form $(e_n,\Im(\log(\Phi(\lambda+
i0)))e_n)_{\calH_+}$
exists for all $\lambda\in\bbR\backslash\calE_n,$ where 
$\calE_n$
has Lebesgue measure zero. Thus one observes,
\begin{align}
&\lim_{\varepsilon\downarrow 0} (e_m,\Im(\log(\Phi_+(\lambda+
i\varepsilon)))e_n)_{\calH_+}
= (e_m,\Im(\log(\Phi_+(\lambda+i0)))e_n)_{\calH_+} \no \\
&=\pi (e_m,\Xi_+(\lambda)e_n)_{\calH_+} \text{ for }
\lambda\in\bbR\backslash\{\calE_m\cup\calE_n\}. \lb{3.33a}
\end{align}
Let $\calE=\cup_{n\in\bbN} \calE_n,$ then $|\calE|=0$ 
($|\cdot|$
denoting the Lebesgue measure on $\bbR$) and hence
\begin{align}
&(f,\Im(\log(\Phi_+(\lambda+i0)))g)_{\calH_{+}} =\pi(f,
\Xi_+(\lambda)g)_{\calH_+} \lb{3.33b} \\
& \text{for } \lambda\in\bbR\backslash\calE
\text{ and } f,g \in \calD=\text{lin.span}\,
\{e_n\in\calH_+\,|\,n\in\bbN\}.
\no
\end{align}
Since $\calD$ is dense in $\calH_+$ and $\Xi_+(\lambda)\in
\calB(\calH_+)$ one infers
 $\Im(\log(\Phi_+(\lambda+i0)))=\pi\Xi_+(\lambda)$ for a.e.
$\lambda \in \bbR,$ completing the proof.
\end{proof}
%%%%%%%%%%%%%%%%%%%%%%%%%%%%%%%%%%%%%%%%%%%%%%%%%%%%%%%%%%%%

Of course \eqref{3.26aa}, \eqref{3.26ab} (and the method of 
proof)
immediately extend to a.e. nontangential limits to the real 
axis.

Assuming Hypothesis~\ref{h3.1} we define
\begin{equation}\lb{3.20}
\xi(\lambda)=
\xi_+(\lambda)-
\xi_-(\lambda)
 \text{ for a.e. }
\lambda\in \bbR
\end{equation}
and call $\xi(\lambda)$ (resp., $\xi_+(\lambda),$ 
$\xi_-(\lambda)$) the
spectral shift function associated with the pair $(H_0,H)$ 
(resp.,
$(H_0,H_+),$ $(H_+,H)$), sometimes also denoted by 
$\xi(\lambda,H_0,H),$
etc., to underscore the dependence on the pair involved.

M.~Krein's basic trace formula
\cite{Kr62} is now obtained as follows.
%%%%%%%%%%%%%%%%%%%%%%%%%%%%%%%%%%%%%%%%%%%%%%%%%%%%%%%%%%
\begin{theorem}\lb{t3.4}
Assume Hypothesis \ref{h3.1}. Then {\rm (}$z \in \bbC 
\backslash
\{\spec (H_0) \cup \spec(H) \}${\rm )}
\begin{equation}
\tr_\calH((H-z)^{-1}-(H_0-z)^{-1})
=-\int_\bbR d\lambda \, \xi(\lambda) (\lambda-z)^{-2}.
\lb{3.21}
\end{equation}
\end{theorem}
%%%%%%%%%%%%%%%%%%%%%%%%%%%%%%%%%%%%%%%%%%%%%%%%%%%%%%%%%%
\begin{proof}
Let $z\in \bbC\backslash \bbR$.
By
\eqref{3.19d} and \eqref{3.19f} we infer
\begin{align}
&\tr_{\calH_+}
(\log (\Phi_+(z)))=
\int_\bbR d\lambda \,\xi_+(\lambda)  (\lambda-z)^{-1} ,
\lb{3.22} \\
&
\tr_{\calH_-}
(\log (\widetilde \Phi_-(z)))=
-\int_\bbR d\lambda \, \xi_-(\lambda) (\lambda-z)^{-1} . 
\lb{3.23}
\end{align}
Adding \eqref{3.13a} and \eqref{3.13b}, differentiating
\eqref{3.22} and \eqref{3.23} with respect to $z$
 proves \eqref {3.21} for $z\in \bbC\backslash \bbR$. The 
result extends to all
 $z\in \bbC \backslash \{\spec (H_0) \cup \spec(H) \}$
 by  continuity
of
$((H-z)^{-1}-(H_0-z)^{-1})$ in $\calB_1(\calH)$-norm.
\end{proof}
%%%%%%%%%%%%%%%%%%%%%%%%%%%%%%%%%%%%%%%%%%%%%%%%%%%%%%

In particular, $\xi(\lambda)$ introduced in 
\eqref{3.20} is
Krein's original spectral shift function (up to 
normalization). As
noted in Section~\ref{s2}, the spectral shift operator
$\Xi_+(\lambda)$ in the particular case $V=V_+,$ and its
relation to Krein's spectral shift function 
$\xi_+(\lambda),$
was first studied by Carey \cite{Ca76} in 1976.
%%%%%%%%%%%%%%%%%%%%%%%%%%%%%%%%%%%%%%%%%%%%%%%%%%%%%%%%%
\begin{remark}\lb{r3.5}
(i) As shown originally by M.~Krein \cite{Kr62}, the 
trace formula
\eqref{3.21} extends to
\begin{equation}\lb{3.24}
\tr(f(H)-f(H_0))
=\int_{\bbR} d\lambda \,\xi(\lambda) f'(\lambda)
\end{equation}
for appropriate functions $f$. This fact has been
studied by numerous authors and we refer, for instance, to
\cite{BW83}, Ch.~19, \cite{BP98}--\cite{BY93},
\cite{Kr83},
\cite{Kr89}, \cite{Pe85},
\cite{Si75},
\cite{SM94},
\cite{Vo87},
\cite{Ya92}, Ch.~8 and the references therein.

\noindent (ii) While we focus here on pairs
of self-adjoint operators $(H_0, H)$, M.~Krein in his 
original paper
\cite{Kr62}
also considered pairs of unitary operators. Given the
conformal equivalence of
$\bbC_+$
 and the open interior of the unit disk $D$ in $\bbC,$
this corresponds precisely to the study of Nevanlinna,
 respectively, Riesz-Herglotz functions in $\bbC_+$, 
respectively,
$D$. Pairs of unitary operators are also studied, for 
instance, in
\cite{BS75},
\cite{BY93},
\cite{MS96},
\cite{Si75},
\cite{SM94}.
The trace formula in the case of non-self-adjoint and 
non-unitary
pairs of operators is a rather challenging problem. The 
interested
reader can get some insight into this matter by consulting
\cite{AP80}, \cite{Kr89}, \cite{La65}, \cite{Ry95}
and the extensive literature cited therein. Similarly, the
case of non-trace
class perturbations and associated trace formulas has
been studied by various authors. We refer, for instance,
to
\cite{Ry96},
where Hilbert-Schmidt perturbations of unitary operators
are treated in depth.

\noindent (iii) While scattering theory for the pair 
$(H_0,H)$ is not
discussed in this
paper, we remark that
$\xi(\lambda)$,  for a.e. $\lambda \in \spec_{ac}(H_0)$
(the absolutely continuous spectrum of $H_0$), is related 
to the scattering
operator at fixed energy $\lambda$ by the Birman-Krein 
formula \cite{BK62},
\begin{equation}\lb{3.25d}
\text{det}_{\calH_\lambda}(S(\lambda, H_0, H))=
e^{-2\pi i \xi(\lambda)}
\text{ for a.e. }
\lambda \in \spec_{ac}(H_0).
\end{equation}
Here
 $S(\lambda, H_0, H)$ denote the fibers in the direct 
integral representation
of the scattering operator
$$
S( H_0, H)=\int_{\spec_{ac}(H_0)}^{\oplus}d\lambda  \,
 S(\lambda, H_0, H)
\text{ in } \calH =
\int_{\spec_{ac}(H_0)}^{\oplus}d\lambda \, 
\calH_\lambda
$$
with respect to the absolutely continuous part 
$H_{0, ac}$ 
of $H_0$.
This fundamental connection, originally due to Birman 
and Krein
\cite{BK62},
is further discussed in
\cite{BW83}, Ch.~19,
\cite{BY93},
\cite{BY93a},
\cite{Ca76},
\cite{Ka78},
\cite{Kr83},
\cite{So93},
\cite{Ya92}, Ch.~8 and the literature cited therein.

\noindent (iv) The standard identity (\cite{GK69}, 
Sect.~IV.3)
\begin{equation}\lb{3.26}
\tr_\calH ((H-z)^{-1}-(H_0-z)^{-1})=
-d \log (\text{det}_\calH (I_\calH+V(H_0-z)^{-1}))/dz
\end{equation}
together with the trace formula \eqref{3.21} yields
the well-known connection between perturbation
 determinants and $\xi(\lambda)$,
 also due to M.~Krein
\cite{Kr62}
\begin{equation}\lb{3.27}
\log (\text{det}_\calH (I_\calH+V(H_0-z)^{-1}))=
\int_\bbR d\lambda \, \xi(\lambda)
 (\lambda -z)^{-1} ,
\end{equation}
\begin{equation}\lb{3.28}
\xi(\lambda)=\lim_{\varepsilon\downarrow 0}\pi^{-1}
\Im (
\log (\text{det}_\calH (I_\calH+V(H_0-(\lambda+ 
i0))^{-1})))
\text{ for a.e. }
\lambda\in \bbR,
\end{equation}
\begin{equation}\lb{3.29}
\tr_\calH(V)=
\int_\bbR d \lambda  \, \xi(\lambda), \quad
\quad
\int_\bbR d \lambda \, \vert \xi(\lambda)\vert \le 
\vert\vert V\vert\vert_1.
\end{equation}
This is discussed in more detail,
for instance, in
\cite{BW83}, Ch.~19,
\cite{BY93},
\cite{Ca76},
\cite{Kr60},
\cite{Kr83},
\cite{KY81},
\cite{Si75}--\cite{So93},
\cite{Ya92}. Relation \eqref{3.28} and the analog of
\eqref{2.2a} for $d(AE_H(\lambda)B)/d\lambda,$ where 
$H=H_0+V,$
$V=B^*A,$ $A,B\in\calB_2(\calH),$ $V=V^*,$ leads to the 
expression
$$
\xi(\lambda)=(-2\pi i)^{-1}\tr_\calH (
\log (I_\calH-2\pi i(I_\calH
-A(H-\lambda-i0)^{-1}B^*)(d(AE_H(\lambda)B^*)/d\lambda)))
$$
for a.e. $\lambda\in\bbR$ (cf., e.g., \cite{BW83}, 
Sects.~3.4.4
and 19.1.4).

\noindent (v) The invariance principle for wave operators of 
the pair
$(H_0, H)$ implies
\begin{equation}\lb{3.30}
\xi(\lambda)=
\xi(\lambda, H_0, H)=
\text{sgn}(\Psi'(\lambda))
\xi(\Psi(\lambda), \Psi(H_0),\Psi(H))
\end{equation}
for a certain class of admissible functions $\Psi$.
 In certain cases (e.g., if $H_0$ and $H$
are bounded from below) this can be used
to define $\xi$ even though $H-H_0=V$ is not of trace class 
as long as
$(\Psi(H)-\Psi(H_0)) \in \calB_1(\calH).$ Prime candidates 
for $\Psi$ in such
cases are semigroup
$(\Psi (\lambda)=e^{-t\lambda}, t>0)$
 and resolvent $(\Psi(\lambda)=(\lambda-z)^{-1},$
$z\in \bbC\backslash \bbR $ or $z<E_0$ for some 
$E_0\in \bbR$) functions.
Pertinent facts in this connection can be found in
\cite{JK78},
\cite{SM94}, and
\cite{Ya92}, Sect.~8.11.

\noindent (vi) For simplicity we chose a single Hilbert 
space formulation
throughout
this section. However, every result immediately extends
to the case where
$K\in \calB_2(\calK, \calH)$,
$J=J^*\in \calB(\calK)$,
$J^2=I_\calK$ and
$\calK$ is another complex separable Hilbert space.
\end{remark}
%%%%%%%%%%%%%%%%%%%%%%%%%%%%%%%%%%%%%%%%%%%%%%%%%%%%%%%%%%%%%%

\begin{remark} \lb{r3.8} Suppose $H_0, H_1,$ and $H_2$ are 
self-adjoint
operators in $\calH$ with $(H_j - H_k)\in\calB_1(\calH)$ 
for all
$j,k\in\{0,1,2\}.$ Denoting by $\xi(\lambda,H_j,H_k)$ the 
Krein spectral
shift function of the pair $(H_j,H_k)$ such that
\begin{equation}
\tr_\calH((H_k-z)^{-1}-(H_j-z)^{-1})=-\int_\bbR d\lambda \,
\xi(\lambda,H_j,H_k)(\lambda-z)^{-2}, \lb{3.45}
\end{equation}
the chain rule,
\begin{equation}
\xi(\lambda,H_0,H_2)=\xi(\lambda,H_0,H_1)+
\xi(\lambda,H_1,H_2)
\text{ a.e.,} \lb{3.46}
\end{equation}
together with
\begin{equation}
\xi(\lambda,H_j,H_k)= -\xi(\lambda,H_k,H_j) \text{ a.e.,} 
\lb{3.47}
\end{equation}
and
\begin{equation}
\xi(\lambda,H_j,H_k)\geq0 \text{ a.e. if } (H_k-H_j)\geq 0, 
\lb{3.48}
\end{equation}
imply the monotonicity property
\begin{equation}
\xi(\lambda,H_0,H_2)\geq \xi(\lambda,H_0,H_1) 
\text{ a.e. if }
H_2\geq H_1. \lb{3.49}
\end{equation}
Here \eqref{3.48} is clear from \eqref{3.16}, \eqref{3.19}.
Equation  \eqref{3.47}
follows from \eqref{3.45}, and \eqref{3.46} is a 
consequence of
\eqref{3.28} observing the facts
\begin{align}
I_\calH+V(H_0-z)^{-1}&=(H_0+V-z)(H_0-z)^{-1}, \quad
V\in\calB(\calH),\no \\
\text{det}_\calH((I_\calH+A)(I_\calH+B))&=
\text{det}_\calH(I_\calH+A)
\text{det}_\calH (I_\calH+B), \quad A,B\in\calB_1(\calH). 
\no
\end{align}
\end{remark}

Given the monotonicity property \eqref{3.49} of Krein's 
spectral
shift function, it is natural to inquire whether or not 
this property
is shared by the spectral shift operator. More precisely, 
one might ask
whether or not
\begin{equation}
\Xi(\lambda,H_0,H_2)\overset{\text{{\bf ?}}}{\geq} 
\Xi(\lambda,H_0,H_1)
\text{ a.e. if } H_2\geq H_1\geq H_0. \lb{3.50}
\end{equation}
The following simple counter example destroys such hopes.

\begin{example} \lb{e3.9}
Let $H_0=0,$ $K\in\calB_2(\calH),$ $J=I_{\calH},$ and hence
$\Phi(z)=I_\calH-K^*Kz^{-1},$ $z\in\bbC\backslash\{0\}.$ Then
\begin{equation}
\lim_{\varepsilon\downarrow 0} \|\log(I_\calH-K^*K(\lambda
+i\varepsilon)^{-1}) - 
\log(I_\calH-K^*K\lambda^{-1})\|_{\calB_1(\calH)}
=0 \lb{3.51}
\end{equation}
for $\lambda\in\bbR\backslash\{\spec(K^*K)\cup\{0\}\}$ and
\begin{equation}
\pi^{-1} \Im(\log(I_\calH -K^*K\lambda^{-1})) = \Xi(\lambda),
 \quad \lambda\in\bbR\backslash\{\spec(K^*K)\cup\{0\}\}.
\lb{3.52}
\end{equation}
Decomposing the self-adjoint operator 
$I_\calH-K^*K\lambda^{-1},$
$\lambda\in\bbR\backslash\{\spec(K^*K)\cup\{0\}\}$ into its 
positive and
negative spectral parts then yields
 $\Xi(\lambda)=\theta(K^*K-\lambda I_\calH),$
 $\lambda\in\bbR\backslash\{\spec(K^*K)\cup\{0\}\},$
where $\theta(\cdot)$
denotes the usual step function (i.e., $\theta(x)=1$ for 
$x>0$
and $\theta(x)=0$ for $x<0$).

Next we choose $\calH=\bbC^2,$ $0<a<b<c<1,$ $ac-b^2\geq 0,$ 
and
\begin{equation}
K_1^*K_1=\begin{pmatrix}1 & b \\ b & 1  \end{pmatrix}, 
\quad
\spec(K_1^*K_1)=\{1-b,1+b\}, \lb{3.54}
\end{equation}
with eigenvectors $2^{-1/2}(1,\pm 1)^t$ associated to the
eigenvalues $1\pm b$ and
\begin{equation}
K_2^*K_2=\begin{pmatrix}1+a & 0 \\ 0 & 1+c  \end{pmatrix}, 
\quad
\spec(K_2^*K_2)=\{1+a,1+c\}, \lb{3.55}
\end{equation}
with eigenvectors $(1,0)^t$ (resp., $(0,1)^t$) associated 
to the
eigenvalue $1+a$ (resp., $1+c$). Finally, choosing
$\lambda\in(1+a,1+b)$ then yields
\begin{align}
\Xi_1(\lambda)&=\theta(K_1^*K_1-\lambda I_{\bbC^2})
=E_{K_1^*K_1}(\{1+b\}), \lb{3.56} \\
\Xi_2(\lambda)&=\theta(K_2^*K_2-\lambda I_{\bbC^2})
=E_{K_2^*K_2}(\{1+c\}), \lb{3.57}
\end{align}
where $\{E_{K^*K}(\lambda)\}_{\lambda\in\bbR}$ denotes 
the family of
orthogonal spectral projections of $K^*K.$
Clearly $H_2=K_2^*K_2\geq K_1^*K_1 = H_1\geq H_0=0$ but
 $\Xi_2(\lambda) \ngeq \Xi_1(\lambda)$ for
$\lambda \in (1+a,1+b)$ since one-dimensional 
self-adjoint projections
cannot satisfy an order relation unless one is a real 
multiple of
the other. (Note, however, that 
$\tr(\Xi_2(\lambda))$$=1+c > 1+b$
$=\tr(\Xi_1(\lambda))$ in accordance with \eqref{3.49}.)
\end{example}

This example shows, in particular, that the chain rule 
\eqref{3.46} 
for $\xi(\lambda)$ does not extend to $\Xi(\lambda).$
%%%%%%%%%%%%%%%%%%%%%%%%%%%%%%%%%%%%%%%%%%%%%%%%%%%%%%%%%%%%
%%%%%%%%%%%%%%%%%%%%%%%%%%%%%%%%%%%%%%%%%%%%%%%%%%%%%%%%%%%

\section{Spectral Averaging: An Operator-Theoretic 
Approach}\lb{s4}

In this section we apply the formalism developed in 
Sections~\ref{s2}
and \ref{s3} to provide an effortless proof
 of spectral averaging and its relation to Krein's spectral 
shift
function as originally  proven by Birman and Solomyak
\cite{BS75}.

For the basic setup of this section we assume the following
hypothesis.
%%%%%%%%%%%%%%%%%%%%%%%%%%%%%%%%%%%%%%%%%%%%%%%%%%%%%%%%%%%%%
\begin{hypothesis}\lb{h4.1}
Let $H_0$ be a self-adjoint operator in
$\calH$ with $\dom (H_0)$, and assume
$\{V(s)\}_{s\in \Omega}\subset \calB_1(\calH)$ to be a 
family of 
self-adjoint trace class operators in $\calH$, where 
$\Omega \subseteq
 \bbR$ denotes an open interval. Moreover, suppose that 
$V(s)$
 is continuously differentiable with respect to $s\in\Omega$
in trace norm.
\end{hypothesis}
%%%%%%%%%%%%%%%%%%%%%%%%%%%%%%%%%%%%%%%%%%%%%%%%%%%%%%%%%%%%
To begin our discussions we temporarily assume that 
$V(s)\geq 0$, that 
is, we suppose 
\begin{equation} \lb{4.2}
V(s)=K(s)K(s)^*, \quad s\in \Omega
\end{equation}
for some $K(s)\in\calB_2(\calH),$ $s\in\Omega.$ Given 
Hypothesis~\ref{h4.1} we define the self-adjoint operator 
$H(s)$ in $\calH$ by
\begin{equation}
H(s)=H_0+V(s), \quad \dom (H(s))=\dom (H_0),
\quad s\in \Omega. \lb{4.3}
\end{equation}

In  analogy to \eqref{3.5}  and \eqref{3.6} we introduce
in $\calH$ $(s\in \Omega,$ $  z\in \bbC \backslash \bbR$),
\begin{equation}\lb{4.7}
\Phi(z,s)
=I_{\calH}+ K(s)^*(H_0-z)^{-1}K(s)
\end{equation}
and hence infer from Lemma \ref{l3.2} that
\begin{equation}\lb{4.9}
\Phi (z,s)^{-1}
=I_\calH-K(s)^*(H(s)-z)^{-1}K(s).
\end{equation}
The following is an elementary but useful result needed 
in the context of
Theorem \ref{t4.3}.
%%%%%%%%%%%%%%%%%%%%%%%%%%%%%%%%%%%%%%%%%%%%%%%%%%%%%%%%%%%
\begin{lemma}\lb{l4.2}
Assume Hypothesis~\ref{h4.1} and \eqref{4.2}. Then
$(s\in \Omega, $ $ z \in \bbC\backslash \bbR),$
\begin{equation}\lb{4.11}
d \tr_{\calH}
(\log (\Phi(z,s)))/ds=
\tr_{\calH}(V'(s)(H(s)-z)^{-1}).
\end{equation}
\end{lemma}
%%%%%%%%%%%%%%%%%%%%%%%%%%%%%%%%%%%%%%%%%%%%%%%%%%%%%%%%%%%
\begin{proof}
 By \eqref{3.15b}, \eqref{4.7},  and \eqref{4.9}
one infers for $z=iy,$ $|y|>0$ sufficiently large,
\begin{align}
&d \tr_{\calH} (\log (\Phi(z,s)))/ds 
=d \tr_{\calH} (\log (I_\calH - V(s)(H_0 -z)^{-1}))/ds \no \\
&=(d/ds)\sum_{j=0}^\infty (-1)^j (j+1)^{-1} 
\tr_{\calH}((V(s)(H_0-z)^{-1})^{j+1}) \no \\
&=\tr_{\calH}(V'(s)((H_0-z)^{-1}\sum_{j=0}^\infty 
(-1)^j (V(s)(H_0 -z)^{-1})^j) \no \\
&=\tr_{\calH}(V'(s)((H(s) -z)^{-1}) \lb{4.1}
\end{align}
by repeated use of \eqref{3.15'} and \eqref{3.12a}. Analytic 
continuation of \eqref{4.1} with respect to 
$z\in\bbC\backslash\bbR$ then proves \eqref{4.11}.
\end{proof}
%%%%%%%%%%%%%%%%%%%%%%%%%%%%%%%%%%%%%%%%%%%%%%%%%%%%%%%%%%%%%

Next, applying Lemma \ref{l3.4} to
$\Phi(z,s)$ in \eqref{4.7} one infers $(s\in \Omega),$
\begin{align}
&\log (\Phi(z,s))=\int d\lambda \, \Xi(\lambda,s)
(\lambda-z)^{-1},
\lb{4.19a} \\
&0\le \Xi (\lambda,s) \le I_\calH, \quad
\Xi (\lambda ,s)\in \calB_1(\calH)
\text{ for a.e. } \lambda \in \bbR, \lb{4.21} \\
&\vert\vert \Xi (\cdot,s)\vert\vert_1 \in L^1(\bbR;d\lambda), 
\no
\end{align}
where $\Xi(\lambda,s)$ is associated with the pair 
$(H_0,H(s)),$ assuming $H(s)\geq H_0,$ $s\in\Omega.$

Our principle result on averaging the spectral measure of
$\{E_{H(s)}(\lambda)\}_{\lambda\in \bbR}$ of $H(s)$ then 
reads 
as follows.

%%%%%%%%%%%%%%%%%%%%%%%%%%%%%%%%%%%%%%%%%%%%%%%%%%%%%
\begin{theorem}\lb{t4.3}
Assume Hypothesis \ref{h4.1} and 
$[s_1,s_2]\subset \Omega$. Let $\xi(\lambda,s)$ be the 
spectral shift function associated with the pair 
$(H_0,H(s))$ {\rm (}cf. \eqref{3.20}{\rm )}, where 
$H(s)$ is defined 
by \eqref{4.3} {\rm (}and we no longer suppose 
$H(s)\geq H_0${\rm )}. Then
\begin{equation} \lb{4.25}
\int_{s1}^{s_2} ds \, ( d (\tr_{\calH} 
(V'(s)E_{H(s)}(\lambda) )))=
(\xi(\lambda, s_2)-\xi(\lambda, s_1))d\lambda.
\end{equation}
\end{theorem}
%%%%%%%%%%%%%%%%%%%%%%%%%%%%%%%%%%%%%%%%%%%%%%%%%%%%%%%%%
\begin{proof}
First we prove \eqref{4.25} in the case $V(s)\geq 0$.
The monotone convergence theorem,
\eqref{4.19a}, and Lemma \ref{l4.2} then yield
$(z\in \bbC\backslash \bbR),$
\begin{align}
&\tr_{\calH} \bigg (\int_\bbR d\lambda \,
((\lambda -\Re (z))^2+
(\Im (z))^2)^{-1}\Im (z)
(\Xi(\lambda, s_2) - \Xi(\lambda, s_1))\bigg )
\no \\
&=\int_\bbR d\lambda \,
((\lambda -\Re (z))^2+
(\Im (z))^2)^{-1}\Im (z)
(\xi(\lambda, s_2)
-
\xi(\lambda, s_1))
\no\\
&=\tr_{\calH} (\Im (\log (\Phi(z,s_2))))
- \tr_{\calH} (\Im (\log (\Phi(z,s_1))))
\no\\
&=\int_{s_1}^{s_2}ds \,
\bigg(\frac{d}{ds} \tr_{\calH}( \Im (\log (\Phi(z,s))))\bigg)
=\int_{s_1}^{s_2}ds
\tr_{\calH} (V'(s)\Im ((H(s)-z)^{-1})). \no
\lb{4.26}
\end{align}
By  the spectral theorem applied to $H(s)$ one obtains
\begin{equation}\lb{4.27}
\Im ((H(s)-z)^{-1})=
\Im (z)
\int_\bbR
d E_{H(s)}(\lambda)
((\lambda -\Re(z))^2 + (\Im (z))^2)^{-1}
\end{equation}
and hence
\begin{align}
&
\int_\bbR d\lambda \,
((\lambda -\Re (z))^2+
(\Im (z))^2)^{-1}
(\xi(\lambda, s_2)
-
\xi(\lambda, s_1))
\no \\
&
=\int_{s_1}^{s_2} ds
\tr_{\calH}\bigg(
\int_\bbR
 d E_{H(s)}(\lambda))
((\lambda -\Re(z))^2 + (\Im (z))^2)^{-1}
V'(s)\bigg). \lb{4.28}
\end{align}
Decomposing the self-adjoint trace class operator $V'(s)$
into its positive and negative parts,
\begin{equation}\lb{4.29}
V'(s)=(V'(s))_+-(V'(s))_-, \quad
0\le (V'(s))_\pm \in \calB_1(\calH),
\end{equation}
the monotone convergence theorem yields
\begin{align}
&
\int_\bbR d\lambda \, 
((\lambda -\Re (z))^2+
(\Im (z))^2)^{-1}
(\xi(\lambda, s_2)
-
\xi(\lambda, s_1))
\no \\
&
=\int_{s_1}^{s_2}ds
\int_\bbR
((\lambda -\Re (z))^2+
(\Im (z))^2)^{-1}
(d(\tr_\calH((V'(s))_+^{1/2} E_{H(s)}
(\lambda) (V'(s))_+^{1/2}) 
\no \\
& \quad -d(\tr_\calH((V'(s))_-^{1/2} E_{H(s)}(\lambda) 
(V'(s))_-^{1/2})))
\no \\
&
=\int_{s_1}^{s_2} ds
\int_\bbR
((\lambda -\Re (z))^2+
(\Im (z))^2)^{-1}
d(\tr_\calH(V'(s)E_{H(s)}(\lambda)))
\no \\
&
=\int_\bbR
((\lambda -\Re (z))^2+
(\Im (z))^2)^{-1}
\int_{s_1}^{s_2} ds (
d(\tr_\calH(V'(s)E_{H(s)}(\lambda)))) \lb{4.30}
\end{align}
using Fubini's theorem in the last step. Thus, by the
uniqueness property of Poisson transforms,
\begin{equation}\lb{4.31}
(\xi(\lambda, s_2)
-\xi(\lambda, s_1))d\lambda =
\int_{s_1}^{s_2} ds \, (
d(\tr_\calH(V'(s)E_{H(s)}(\lambda)))).
\end{equation}

In the case of arbitrary $V(s)\in\calB_1(\calH)$ (not 
necessarily 
nonnegative), we argue as follows. Define (cf. \eqref{4.29})
\begin{equation}
W=\bigg(\int_{s_1}^{s_2} ds\, (V'(s))_-\bigg) -V(s_1), 
\lb{4.32}
\end{equation}
then $W\in\calB_1(\calH)$ and $V(s)+W\geq 0$ for all 
$s\in [s_1,s_2].$ Equation \eqref{4.25} now follows from 
the chain 
rule \eqref{3.46} for spectral shift functions
\begin{equation}
\xi(\lambda,H_0, H(s))= \xi(\lambda,H_0, H_0-W)+
\xi(\lambda,H_0-W, H(s)).
\end{equation}
Indeed, $(V(s)+W)'=V'(s),$ and $\xi(\lambda,H_0,H_0-W)$ 
is independent of $s$ and hence drops out on the 
right-hand side of \eqref{4.25}. Moreover, the 
pair $(H_0-W,H(s))$ only involves the nonnegative perturbation 
$V(s)+W$ of $H_0-W$ in $H(s)=H_0-W+(V(s)+W)$ so 
that Lemma~\ref{l4.2} 
becomes applicable as in \eqref{4.31} in the first part 
of our proof.
\end{proof}
%%%%%%%%%%%%%%%%%%%%%%%%%%%%%%%%%%%%%%%%%%%
\begin{remark}\lb{r4.4}
 (i) In the special case of averaging over the boundary 
condition
 parameter for half-line Sturm-Liouville operators 
(effectively a
rank-one resolvent perturbation problem), Theorem \ref{t4.3} 
has
 first been derived by Javrjan
\cite{Ja66},
\cite{Ja71}. The case of rank-one perturbations was 
recently treated
in detail by Simon \cite{Si95}.
The general case of trace class perturbations is due to 
Birman and
 Solomyak
\cite{BS75} using an approach of Stieltjes' double operator 
integrals.
Birman and Solomyak treat the case $V(s)=sV,$ 
$V\in\calB_1(\calH),$ $s\in [0,1].$ As explained in
\cite{Bi94},
\cite{BS75},
and
\cite{BY93},
 the authors were interested in a real analysis
approach to the
spectral shift function in contrast to M.~Krein's
complex analytic treatment. In this way the  local
 integrability of $\xi$  could not be obtained directly
(only its property of  being  a generalized function 
was obtained)
although it follows of course from the uniqueness  $\xi$
(up to additive constants). While our operator  theoretic 
approach
 is intrinsically
 complex-analytic, and hence very much in M.~Krein's spirit,
 it leads to a natural proof of the absolute
 continuity of the (signed) measure on the right-hand
side of \eqref{4.25}. A short proof of \eqref{4.25}
(assuming $V'(s)\ge 0$) has recently been given by Simon
\cite{Si98}.

\noindent (ii) We note that variants  of \eqref{4.25} in  
the context of
one-dimensional
 Sturm-Liouville operators (i.e.,
variants of Javrjan's results in
\cite{Ja66},
\cite{Ja71})
have been repeatedly rediscovered by several authors. In 
particular,
the absolute continuity of averaged spectral measures (with 
respect to
 boundary condition parameters or coupling constants of 
rank-one perturbations)
has been used to prove localization properties
 of one-dimensional random Schr\"odinger operators 
(see, e.g.,
\cite{BS98}, \cite{Ca83},
\cite{Ca84}, \cite{CL90}, Ch.~VIII, \cite{CH94},
\cite{CHM96}, \cite{KS98}--\cite{KS87},
\cite{PF92}, Ch.~V, \cite{Si85}, \cite{Si95}).

\noindent (iii) We emphasize that Theorem~\ref{t4.3} 
applies 
to unbounded
operators (and
hence to random Schr\"odinger operators bounded from below) 
as long as
appropriate relative trace class conditions (either with 
respect to resolvent
or semigroup perturbations) are satisfied.

\noindent (iv) In the special case
$V'(s)\geq 0,$ the  measure 
$$
d(\tr_\calH (V'(s)E_{H(s)}(\lambda))) =
d(\tr_\calH (V'(s)^{1/2}E_{H(s)}(\lambda)V'(s)^{1/2}))
$$
 in \eqref{4.25}
represents a positive measure.

\noindent (v) The result \eqref{4.25} is not restricted to 
a one-dimensional parameter space $s\in[s_1,s_2].$ In fact, 
if $\gamma({\bf s_1},{\bf s_2})$ denotes an oriented 
piecewise $C^1$-path connecting 
${\bf s_1}\in\bbR^n$ and ${\bf s_2}\in\bbR^n,$ one obtains 
analogously,
\begin{equation}
\int_{\gamma({\bf s_1},{\bf s_2})} d{\bf s}\cdot 
(d(\tr_\calH(({\bf \nabla}V)({\bf s})
E_{H(s)}(\lambda)))) = (\xi(\lambda,{\bf s_2}) - 
\xi(\lambda,{\bf s_1}))d\lambda. \lb{4.32a}
\end{equation}
We omit further details.
\end{remark}
%%%%%%%%%%%%%%%%%%%%%%%%%%%%%%%%%%%%%%%%%%%%%%%%%%%%%%%%%%%%

In the special case of a sign-definite perturbation of $H_0$ 
of the form
$sKK^*,$ one can in fact prove an
operator-valued averaging formula as follows.
%%%%%%%%%%%%%%%%%%%%%%%%%%%%%%%%%%%%%%%%%%%%%%%%%%%%%%%%%%%%%
\begin{theorem} \lb{t4.5}
Assume Hypothesis~\ref{h3.1} and $J=I_\calH.$ Then
\begin{equation}
\int_0^1 ds \, d(K^*E_{H_0 + sKK^*}(\lambda)K)=
\Xi(\lambda)d\lambda, \lb{4.33}
\end{equation}
where $\Xi(\cdot)$ is the spectral shift operator 
associated with
\begin{equation}
\Phi(z)=I_\calH + K^*(H_0 - z)^{-1}K, 
\quad z\in\bbC\backslash\bbR,
\lb{4.34}
\end{equation}
that is,
\begin{align}
&\log(\Phi(z))= \int_\bbR d\lambda \, \Xi(\lambda) 
(\lambda - z)^{-1},
\quad
z\in\bbC\backslash\bbR, \lb{4.35} \\
& 0\leq \Xi(\lambda) \in \calB_1 (\calH) 
\text{ for a.e. }\lambda\in\bbR,
\quad \|\Xi(\cdot)\|_1\in L^1(\bbR;d\lambda). \lb{4.36}
\end{align}
\end{theorem}
%%%%%%%%%%%%%%%%%%%%%%%%%%%%%%%%%%%%%%%%%%%%%%%%%%%%%%%%%%
\begin{proof}
An explicit computation shows
\begin{align}
&(\lambda - \Phi(z))^{-1} = -(1-\lambda)^{-1} \no \\
&\times
(I_\calH
-(1-\lambda)^{-1}K^*(H_0+(1-\lambda)^{-1}KK^* - z)^{-1}K) 
\in\calB(\calH)
 \lb{4.37}
\end{align}
for all $\lambda <0.$ Since $\log(\Phi(z))=\ln(\Phi(z))$ for
$z\in\bbC_+$ as a result of analytic continuation, one 
obtains
\begin{align}
\log(\Phi(z)) &= \int_\bbR d\lambda \, \Xi(\lambda) 
(\lambda -z)^{-1}
\no \\
=\ln(\Phi(z)) &= \int_{-\infty}^0 d\lambda \, 
((\lambda -\Phi(z))^{-1} -
\lambda(1+\lambda^2)^{-1}I_\calH) \no \\
&=\int_{-\infty}^0
d\lambda \, (1-\lambda)^{-2} K^*
(H_0+(1-\lambda)^{-1}KK^*-z)^{-1}K \no \\
&= \int_0^1 ds \int_\bbR
d(K^*E_{H_0+sKK^*}(\lambda)K)(\lambda-z)^{-1} \no \\
&=\int_\bbR (\lambda-z)^{-1}\int_0^1
ds \, d(K^*E_{H_0+sKK^*}(\lambda)K) \lb{4.38}
\end{align}
proving \eqref{4.33}. (Here the interchange of the 
$\lambda$ and $s$
integrals follows from Fubini's theorem considering 
\eqref{4.38}
in the weak sense.)
\end{proof}
%%%%%%%%%%%%%%%%%%%%%%%%%%%%%%%%%%%%%%%%%%%%%%%%%%%%%%%%%%

As a consequence of Theorem~\ref{t4.5} one obtains
\begin{equation}
\int_{s_1}^{s_2} ds\, d(K^* E_{H_0+sKK^*}(\lambda)K) =
\Xi(\lambda,s_2) - \Xi(\lambda,s_1), \lb{4.39}
\end{equation}
where $\Xi(\lambda,s)$ is the spectral shift operator 
associated with
$\Phi(z,s)=I_\calH + sK^*(H_0-z)^{-1}K,$ $s\in [s_1,s_2].$

%%%%%%%%%%%%%%%%%%%%%%%%%%%%%%%%%%%%%%%%%%%%%%%%%%%%%%%%%%%
\vspace*{3mm}
\noindent {\bf Acknowledgments.}
F.~G. would like to thank P.~Exner for his kind invitation 
to take part
in this conference, and all organizers for providing a
most stimulating meeting. He also thanks H.~Holden for 
the generous
invitation
to spend seven weeks in the summer of 1998 at the 
Norwegian University of
Science and Technology, NTNU, Trondheim, during the final 
stages of
this work. The extraordinary hospitality at the 
Department of
Mathematical Sciences at NTNU and financial support by 
the
Research Council Norway, Grant No.~107510/410 are 
gratefully
acknowledged.

K.~A.~M. gratefully acknowledges financial support by
 the  University of Missouri Research Board, Grant 
No.~98--119.

S.~N.~N. thanks H.~Holden for making it possible to spend 
a week in July
of 1998 at the Norwegian University of Science and 
Technology, NTNU,
Trondheim.

\end{document}